\documentstyle{amsppt}
\magnification 1200
\def\today
{\ifcase\month\or
     January\or February\or March\or April\or May\or June\or
     July\or August\or September\or October\or November\or December\fi
     \space\number\day, \number\year}
\magnification 1200
\input pictex
\UseAMSsymbols
\hsize 5.5 true in
\vsize 8.5 true in
\parskip=\medskipamount
\NoBlackBoxes

\def\mathbb{\Bbb}

\def\mathcal{\Cal}

\def\supp{\text{\rm supp\,}}

\def\mod{\text{\rm mod\,}}

\def\ve{\varepsilon}
\def\vp{\varphi}
\def\arrowk{^\to{\kern -6pt\topsmash k}}
\def\arrowK{^{^\to}{\kern -9pt\topsmash K}}
\def\arrowr{^\to{\kern-6pt\topsmash r}}
\def\bark{\bar{\kern-0pt\topsmash k}}
\def\arrowvp{^\to{\kern -8pt\topsmash\vp}}
\def\arrowf{^{^\to}{\kern -8pt f}}
\def\arrowg{^{^\to}{\kern -8pt g}}
\def\arrowu{^{^\to}a{\kern-8pt u}}
\def\arrowt{^{^\to}{\kern -6pt t}}
\def\arrowe{^{^\to}{\kern -6pt e}}
\def\tk{\tilde{\kern 1 pt\topsmash k}}
\def\barm{\bar{\kern-.2pt\bar m}}
\def\barN{\bar{\kern-1pt\bar N}}
\def\barA{\, \bar{\kern-3pt \bar A}}

\def\mathbb{\Bbb}

\def\nint{\not\!\!\int}
\def\dist{\text{\rm dist\,}}

\TagsOnRight
\NoRunningHeads

\def\supp{\text{\rm supp\,}}

\def\mod{\text{\rm mod\,}}

\def\ve{\varepsilon}
\def\vp{\varphi}
\hsize = 6.2true in

\vsize=8.2 true in
\TagsOnRight
\NoRunningHeads

\document
\topmatter
\title
Bounds on Oscillatory Integral Operators Based on Multilinear Estimates
\endtitle
\author
J.~Bourgain and L.~Guth
\endauthor
\address
Institute for Advanced Study, Princeton, NJ 08540
\endaddress
\email
bourgain\@ias.edu, lguth\@ias.edu
\endemail
\endtopmatter

\centerline
{\today}

\noindent
{\bf \S1. Summary} 

Let $S\subset \mathbb R^n$ be a smooth, compact hyper-surface with positive definite second fundamental form.
Let $\sigma$ be its surface measure.

We prove the following result with respect to the Fourier restriction/extension problem.

\proclaim
{Theorem 1} Assume the exponent $p$ satisfies
$$
\left\{
\aligned
p> 2\frac{4n+3}{4n-3} &\text { if } n\equiv 0 \ (\text{mod } 3)\\
p> \frac{2n+1}{n-1}  &\text { if } n\equiv 1 \ (\text{mod } 3)\\
p > \frac {4(n+1)}{2n-1} &\text { if } n\equiv 2 \ (\text{mod } 3).
\endaligned
\right.
\tag 1.1
$$
Then the inequality
$$
\Vert\hat\mu\Vert_p \lesssim C_p \Big\Vert\frac{d\mu}{d\sigma}\Big\Vert_\infty\tag 1.2
$$
holds for measures $\mu \ll \sigma $ such that $\frac {d\mu}{d\sigma} \in L^\infty (S, \sigma)$.
\endproclaim

See \S3. 
For $n=3$ (resp. $n=4$), the exponent in (1.2) is $\frac{10}3$ (resp. 3) and coincides with the condition
$p\geq \frac{2(n+2)}n$ resulting from the bilinear $L^2$-approach in [T1].
For $n\geq 5$, the result is new.

Recall that, according to the restriction conjecture, due to E.~Stein, cf. [St1], (1.1) should remain valid for all $p>\frac {2n}{n-1}$.

We also  point out that if $S$ is the $(n-1)$-sphere or paraboloid, then (1.2) may be strengthen to
$$
\Vert\hat \mu\Vert_p \leq C_p \Big\Vert\frac {d\mu}{d\sigma}\Big\Vert_p\tag 1.2$'$
$$
for $p$ satisfying (1.1)
(the argument combines Theorem 1, the Maurey-Nikishin factorization theorem and invariance considerations, the usual way;
cf. [B1]).

The main ingredient in our approach is the multilinear theory developed in [BCT] that we will recall in \S5.
In \S2 we treat the case $n=3$ to explain the method in its simplest form.
In \S4, the analysis is refined further and combined with T.~Wolff's Kakeya maximal function estimate [Wo1] to establish
(1.1) for $n=3$ under the condition
$$
p > 3 \frac 3{10}.\tag 1.3
$$

Thus we have the following small improvement of the $p>\frac {10} 3$ result in 3D.

\proclaim
{Theorem 2}
For $n=3$ and $S$ as above, we have
$$
\Vert\hat\mu\Vert_p \leq C_p \big\Vert\frac{d\mu}{d\sigma}\Big\Vert_\infty \text { for } p> 3\frac 3{10}\tag 1.4
$$
assuming $\mu\ll \sigma$ and $\frac {d\mu}{d\sigma} \in L^\infty (S, d\sigma)$.
\endproclaim

By using `$\ve$-removal lemmas', Theorems 1 and 2 may be derived from a weaker `local' version, more precisely

\proclaim 
{Theorem 1$'$} Let $n\geq 3$ and $S$ as above.

Denote
$$
Q_R^{(p)} =\max\Vert\hat\mu\Vert_{L^p(B_R)}
$$
where the maximum is taken over all measures $\mu\ll\sigma$ on $S$ such that $\Vert\frac{d\mu}{d\sigma}\Vert_\infty \leq 1$.
Then, for all $\ve>0$
$$
Q_R^{(p)} \ll R^\ve\tag 1.5
$$
provided $p$ satisfies (1.1).
\endproclaim
\noindent
and

\proclaim
{Theorem 2$'$} Same statement for $n=3$ and $p\geq 3\frac 3{10}$.
\endproclaim

The use of such $\ve$-removal lemmas is by now standard (cf. [T2]), but we will include an argument for completeness sake in the Appendix,
since we process here $L^\infty -L^p$ inequalities rather than $L^p-L^p$ inequalities, as in [T2].

The technique used applies also in the variable coefficient (H\"ormander) setting.
Thus we consider oscillatory integral operators
$$
(T_\lambda f)(x) =\int e^{i\lambda\psi(x, y)} f(y) dy\qquad (\Vert f\Vert_\infty \leq 1)\tag 1.6
$$
with real analytic phase function
$$
\psi(x, y) =x_1y_1 +\cdots+ x_{n-1} y_{n-1} +x_n\langle A y, y\rangle +O(|x| \, |y|^3)+ O (|x|^2|y|^2)
\tag 1.7
$$
and $A$ non-degenerate.

($x\in\Bbb R^n, y\in \Bbb R^{n-1}$ are restricted to a neighborhood of $0$.)

Our concern is then in which range of $p$, a bound
$$
\Vert T_\lambda f\Vert_p \leq c\lambda^{-\frac np} \tag 1.8
$$
holds.
Recall Stein's result [St2]
$$
\Vert T_\lambda f\Vert_p \leq c\Vert f\Vert_2 \text { for }  p\geq \frac{2(n+1)}{n-1}.\tag 1.9
$$
Also, for $n$ {\it odd}, there are examples showing that, replacing $\Vert f\Vert_2$ by $\Vert f\Vert_\infty$, an inequality (1.8) may
only hold for $p\geq \frac{2(n+1)}{n-1}$ (see [B2]).

Lee observed in [L] that Stein's estimate may be improved if we make the additional hypothesis that 
$A$ in (1.7) is positive (or negative) definite.  He extended the bilinear approach from [T1] to the variable coefficient
setting.  In particular, he proved that (1.8) holds (up to a factor $\lambda^{\epsilon}$) if $p \geq \frac{2(n+2)}{n}$.  We will prove in \S5 that (1.8) 
holds under the condition (1.1). Thus we have

\proclaim
{Theorem 3}
Let $T_\lambda$ be as above with $A$ positive or negative definite in (1.7). Then
$$
\Vert T_\lambda f\Vert_p \leq C_p \lambda^{-\frac np}\Vert f\Vert_\infty\tag 1.10
$$
holds for $p$ satisfying (1.1).
\endproclaim

If $n =3 $ or $4$, Theorem 3 agrees with the results of [L], and for $n \ge 5$, it is new.

For $n$ {\underbar {even}}, there is the following statement (with only the non-degeneracy assumption on $A$).

\proclaim
{Theorem 4} Let $n$ be even and $T_\lambda$ as above, assuming in (1.7) that $A$ is non-degenerate.
Then
$$
\Vert T_\lambda f\Vert_p \leq C_p \lambda^{-\frac np} \Vert f\Vert_\infty \text { for } \ p>\frac {2(n+2)}{n}\tag 1.11
$$
\endproclaim
\noindent
(apart from the endpoint, the condition on $p$ in Theorem 4 was already previously observed to be best possible, cf. [B2].)

It turns out, rather surprisingly, that for $n=3$ the exponent $\frac {10}3 $ in Theorem 3 is also optimal.
In \S6, we describe a specific example (with $A$ elliptic), making the comparison with the hyperbolic case, and explaining the
role of the Kakeya compression phenomenon.
For $n=3$, in both elliptic and hyperbolic cases, there may be a curved Kakeya compression in a 2-dimensional set at
the coarse scale $\frac 1{\sqrt\lambda}$, but the local behaviour of the oscillatory integrals is different.

The proof of Theorems 3 and 4 is based on an application of Theorem 6.2 from [BCT], but we need a version without the extra
$\lambda^\ve$-factors.
Hence, we proceed to `$\ve$-removal' at the multilinear stage (see Appendix), which also provides an alternative strategy to derive Theorem 1 directly,
without passing through Theorem 1$'$ (let us point out that this $\ve$-removal argument applies only
to our particular application of [BCT], Theorem 6.2, see \S5.)

Returning to curved Kakeya compression, it is shown that a curved Kakeya set in even dimension $n$ has Minkowski dimension at least $\frac n2+1$ (see
\S6).
This statement was known to be optimal (see [B2]).

Details are given in \S7 for $n=4$, where it is shown how to derive this property from multi-linear Kakeya-type results.
This strategy may be seen as the essence of our paper and is basically repeated to obtain the oscillatory integral bounds cited above.

Returning to Theorem 3, we should point out the application to the Bochner-Riesz multilinear problem.
Recall that the Bochner-Riesz multiplier $S_\delta$ is defined by $(S_\delta f)^\wedge (\xi)= (1-|\xi|^2)^\delta_+ \hat f(\xi)$.
Equivalently $S_\delta f= f* K_\delta$, where $K_\delta$ has the asymptotic
$$
K_\delta(x) \sim e^{e^{\pm 2\pi i|x|}}/ |x|^{\frac{n+1}{2} +\delta}.\tag 1.12
$$
The problem is then to obtain the optimal condition on $\delta\geq 0$ to satisfy
$$
\Vert S_\delta f\Vert_{L^p(\Bbb R^n)}\leq C\Vert f\Vert_{L^p(\Bbb R^n)}.\tag 1.13
$$

C.~Fefferman's proof of the ball-multiplier conjecture implies that certainly $\delta>0$ for $p\not=2$ (note that the problem is self-dual).
In view of (1.12), the condition
$$
\delta>\max \Big(0, \Big|\frac 12 -\frac 1p\Big| n -\frac 12\Big)\tag 1.14
$$
is clearly necessary.
It is conjectured that (1.14) also suffices for (1.13) to hold and  this was proven for $n=2$ in [C-S] and, independently, in [Hor].

In fact, H\"ormander's approach consists in reducing the study of convolution by $K_\delta$ to some specific oscillatory integral operator $T_\lambda$,
of the type considered above (note that regarding dimension, the $\Bbb R^d -\Bbb R^d$ problem is replaced by an $\Bbb R^{d-1}-\Bbb R^d$ problem
in this reduction).  As a corollary of our Theorem 3 together with the standard factorization and rotational invariance considerations (already
mentioned above), we obtain (cf. [B2] for details).

\proclaim
{Theorem 5} Let $n\geq 3$. Then the Bochner-Riesz conjecture holds providing $\max (p, p')$ satisfies (1.1).
\endproclaim

On the geometric side, the Kakeya-type maximal function underlying the Bochner-Riesz operators (sometimes called `Nikodym maximal function) involves
also averaging over straight line segments and, for $n=3$, T.~Wolff's $\frac 52$-inequality is again known to hold (see [Wo1]).
Thus in principle, one could expect the proof of Theorem 2 to carry over and lead to the validity of the Bochner-Riesz conjecture for
$\max (p, p')\geq 3 \frac 3{10}$, if $n=3$.
We do not pursue the details of this matter here.
In fact, it is well-possible that the exponent $3\frac 3{10}$ from Theorem 2 may be improved further, by reorganizing and refining the method.
No serious attempt was given to do so, as our primary goal is to show how to obtain  some progress over the present results, keeping the arguments
as simple as possible.

Finally, let us cite [T3] as a survey work on the problems discussed in this paper and where the reader will find many background material and
references.

{\bf Acknowledgement.}  The first author was partially supported by NSF grant DMS-0808042.
The second author was supported by NSF grant DMS-0635607 and the Monell Foundation.

The authors are most grateful to the referee for a careful reading and commenting on an earlier manuscript, that led to many improvements.

\bigskip

\noindent
{\bf \S2. An Approach to the Restriction Problem in 3D}

{(alternative proof of the $L^{10/3}$-bound)}
\bigskip

\noindent
{\bf 1.} Consider the oscillatory integral operator
$$
Tf(x) =\int e^{i\phi(x, y)} f(y) dy
\qquad (|f|\leq 1)
$$
where $y\in \Omega$ is a neighborhood of $ 0\in\Bbb R^2$ and $x\in\mathbb R^3\cap [|x|<R]$,
$$
\phi(x, y) =x_1y_1+x_2y_2+ x_3\phi_1 (y)\tag 1.1
$$
with $\phi_1(y) =y^2_1+ y_2^2$ (paraboloid), or more generally
$$
\phi_1(y) =\langle Ay, y\rangle +O(|y|^3) \quad \text {($A$ = positive definite)}
\tag 1.2
$$
(we will comment on the indefinite case at the end of this section).

The purpose of this section is to explain in a simple case how the multi-linear theory from [BCT]
can be exploited to produce results in the usual restriction problem.

Given a phase function $\phi$ as above, we introduce at a given point $y\in\Omega$ the vector
$$
Z= Z(y)=\partial_{y_1} (\nabla_x\phi)\wedge \partial_{y_2} (\nabla_x\phi)
=(-\partial_1 \phi_1(y), -\partial_2 \phi_1 (y), 1).\tag 1.3
$$
For simplicity, we carry the discussion for the case of the paraboloid, thus 
$$
\phi_1(y) =y^2_1+ y_2^2.
$$
In this case, the transversality condition of $\{Z(y^{(i)}), i= 1, 2, 3\}$, where
$y^{(i)}$ is restricted to some small disc $\Omega_i \subset\Omega$ (as needed for the trilinear $L^3$-bound
from [BCT]) amounts to non-collinearity of $\Omega_1, \Omega_2, \Omega_3$.

Discussion of the general situation (1.2) would require to introduce the Gauss map associated to the surface
$$
(y_1, y_2)\mapsto \big(y_1, y_2, \phi_1 (y)\big).
$$
(see \S3.)
\medskip

\noindent
{\bf 2.} \ Fix $K$ (a large parameter).

Partition $\Omega=\bigcup\Omega_\alpha, \Omega_\alpha$ balls of size $\frac 1K$; $y_\alpha\in\Omega_\alpha$.
There are $\sim K^2$ values of $\alpha$.
Write
$$
\align
Tf(x) &=\sum_\alpha e^{i\phi(x, y_\alpha)} \Big[\int_{\Omega_\alpha} e^{i[\phi(x, y)-\phi(x, y_\alpha)]} f(y) dy\Big]
\\
&=\sum_\alpha e^{i\phi(x, y_\alpha)} (T_\alpha f)(x).  \tag 2.1
\endalign
$$

Note that
$$
|\nabla_x[\phi(x, y)-\phi(x, y_\alpha)]|\leq \frac 1K \text{ for $y\in\Omega_\alpha$}.
$$
Take a smooth rapidly decaying bumpfunction $\eta$ s.t. $\hat\eta(\omega) =1$ on
$ [\omega\in\mathbb R^3; |\omega|~\leq~1]$. Let $\eta_K(x) =\frac 1{K^3}\eta \big(\frac xK\big)$ satisfying $\hat\eta_K(\omega) =1$ for 
$|\omega|< 1/K$.

Thus
$$
T_\alpha f=T_\alpha f*\eta_K
$$
and
$$
|T_\alpha f(x) |\leq \int |T_\alpha f(z)| \ |\eta_K(x-z)|dz.
$$
Restrict $x$ to a ball $B(a, K)\subset\mathbb R^3$. Set $a=0$.

For $x\in B(0, K)$
$$
|T_\alpha f(x)|\leq \int |T_\alpha f(z)| \zeta_K(z)dz =c_\alpha\tag 2.2
$$
where
$$
\zeta(x) = \max_{|x-x'|\leq 1} |\eta(x')|.
$$
\bigskip

\noindent
{\bf 3.}  Denote $c_* =\max c_\alpha =c_{\alpha_*}$.
Let $K_1\ll K$ be a second large parameter.  We distinguish several possibilities.

\noindent
(3.1) {\bf Non-coplanar interaction}. 

There are $\alpha, \beta, \gamma$ such that $c_\alpha, c_\beta, c_\gamma>K^{-4} c_*$ and
$$
|y_\alpha -y_\beta|\geq |y_\alpha -y_\gamma| \geq \dist \big
(y_\gamma, \underbrace  {y_\alpha+\mathbb R(y_\beta- y_\alpha)}_{\equiv\ell(y_\alpha, y_\beta)}
\big) > 10^3\frac 1K\tag 3.1$'$
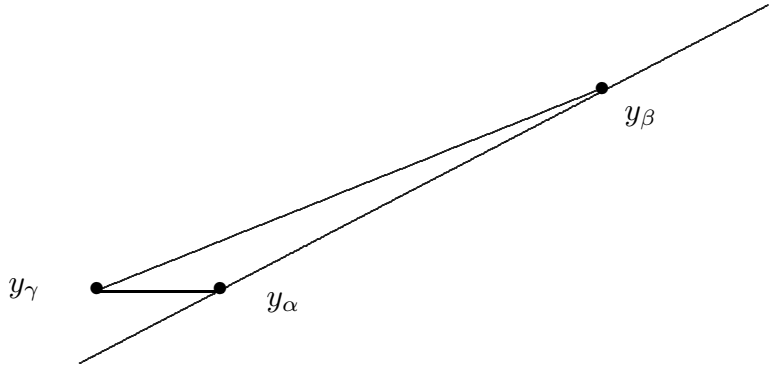
$$
\vskip.3 true in
%%User: gustafss@miami.math.ias.edu (Elly Gustafsson)
\font\thinlinefont=cmr5
$$
\hbox{\beginpicture
\setcoordinatesystem units <1.00000cm,1.00000cm>
%\setshadesymbol ({\thinlinefont .})
\setlinear
%
% Fig POLYLINE object
%
\linethickness= 0.500pt
\setplotsymbol ({\thinlinefont .})
{\plot  4.445 21.431 12.065 25.400 / }%
%
% Fig POLYLINE object
%
\linethickness= 0.500pt
\setplotsymbol ({\thinlinefont .})
{\putrule from  4.604 22.225 to  6.032 22.225 }%
\linethickness= 0.500pt
\setplotsymbol ({\thinlinefont .})
{\plot  4.604 22.225 10.160 24.448 / }%
\put{$\bullet$} [lB] at  4.539 22.170
\put{$\bullet$} [lB] at 10.139 24.390
%
% Fig TEXT object
%
\put{$\bullet$} [lB] at  5.909 22.170
%
% Fig TEXT object
%
\put{$y_\gamma$} [lB] at  3.651 22.225
%
% Fig TEXT object
%
\put{$y_\alpha$} [lB] at  6.509 22.066
\put{$y_\beta$} [lB] at 10.478 24.130
\linethickness=0pt
\putrectangle corners at  3.651 25.425 and 12.090 21.406
\endpicture}
$$

\vskip .3 true in

In this situation we use the trilinear theory from [BCT].

\bigskip
(3.2) \ {\bf Non-transverse interaction}. 

If $|y_\alpha - y_{\alpha_*}| >\frac 1{K_1}$, then $c_\alpha\leq K^{-4} c_*$.
Here we use rescaling (cf. [T-V-V]).

(3.3) \ {\bf Transverse coplanar interaction}.

There is $\alpha_{**}$ with $c_{\alpha_{**}} >K^{-4} c_*, |y_{\alpha_*}-y_{\alpha_{**}}|>\frac 1{K_1}$.

Assuming (3.1) fails, it follows that moreover
$$
c_\alpha\leq K^{-4} c_* \text { if } \dist \big(y_\alpha, \ell(y_{\alpha_*}, y_{\alpha_{**}})\big)> 10^3 \frac {K_1}K.
$$

In this case we rely on the by now standard square function estimates going back to A.~Cordoba's work [C].

\noindent
{\bf 4.}  Assume (3.1)

For $x\in B(0, K)$, by (2.2), (1.1)
$$
|Tf(x)|\leq \sum_\alpha c_\alpha < K^2 c_* <K^6 (c_\alpha c_\beta c_\gamma)^{\frac 13}.
$$
Hence, for $q\geq 3$
$$
\align
|Tf(x)|^q \leq |Tf(x)|^3 &\leq K^{18}\int |T_\alpha f|(z_1) |T_\beta f|(z_2)|T_\gamma f|(z_3) \, \zeta_K(z_1) \zeta_K(z_2)\zeta_K(z_3) d{z_1} d{z_2} d{z_3}\\
&\leq K^{18} \sum_{\alpha, \beta, \gamma \, (3.1')} \int|T_\alpha f|(x-z_1)|T_\beta f|(x-z_2) |T_\gamma f |(x-z_3) \zeta_K(z_1)\zeta_K(z_2)\zeta_K(z_3)
\endalign
$$

The corresponding contribution is estimated using the trilinear bound from [BCT]
$$
\int_{B_R} |T_\alpha f| (x-z_1)|T_\beta f|(x-z_2)| T_\gamma f|(x-z_3) dx< R^\ve. C(K)< R^{2\ve}\tag 4.1
$$

\noindent
{\bf 5.}  Assume (3.2). For $x\in B(0, K)$, estimate
$$
\align
|Tf(x)| &\leq 10\max_\tau \Big|\int_{\tilde \Omega_\tau} e^{i\phi(x,y)} f(y) dy\Big|+\sum_{|y_\alpha -y_{\alpha_*}|>\frac 1{K_1}} c_\alpha\\
&\leq 10. \max_\tau|\tilde T_\tau f(x) | + K^{-2} c_*
\tag  5.1
\endalign
$$
where $\Omega =\bigcup \tilde\Omega_\tau$ is a partition of $\Omega$ in balls of size $\frac 1{K_1}$.
\bigskip

Thus (5.1) implies for $x\in B(0, K)$
$$
|Tf(x)|^q\leq C \sum_{\tau=1}^{\sim K_1^2} |\tilde T_\tau f|^q (x)+C K^{-2q}\sum^{\sim K^2}_{\alpha=1} \int|T_\alpha f|^q (x-z)\zeta_K(z)dz.\tag 5.2
$$

The corresponding contribution is at most
$$
C \sum_\tau\int_{B_R} |\tilde T_\tau f|^q +C K^{-2q} \sum_\alpha \int_{B_R} |T_\alpha f|^q.\tag 5.3
$$
At this point, we use the (parabolic) rescaling
$$
\align
&\Big| \int_{|y-\bar y|<\rho} e^{i\phi(x,y)} f(y) dy\Big|= \\
&{}\\
&y=\bar y+y'\\
&\\
&\Big|\int_{|y'|<\rho} e^{i[(x_1+ 2\bar y_1 x_3) y_1'+(x_2+ 2\bar y_2 x_3) y_2'
+ x_3 |y'|^2]} f(\bar y+y')dy'\Big| = (5.4)
\endalign
$$
and
$$
\Vert(5.4)\Vert_{L^q(B_R)} \leq C\rho^2 \rho^{-\frac 4q}Q_{\rho R}\tag 5.5
$$
where we define
$$
Q_R =\max_{|f|\leq 1} \Vert Tf\Vert_{L^q(B_R)}. \tag 5.6
$$
Substituting (5.5) in (5.3) gives the contribution $(\rho=\frac 1{K_1}$ and $\rho=\frac 1K)$
$$
CK_1^2.K_1^{-2q+4} Q^q_{R/K_1}+ CK^{-2q}. K^{2}. K^{-2q+4} Q^q_{R/K}
$$
and hence for the $L^q$-norm
$$
< CK_1^{-2(1-\frac 3q)} Q_{R/K_1} + CK^{-4+\frac 6\ve} Q_{R/K} .\tag 5.7
$$
\bigskip

\noindent
{\bf 6.}  Assume (3.3).  Thus, denoting $\ell=\ell(y_{\alpha_*}, y_{\alpha_{**}})$, for $x\in B(a, R)$
$$
\align
\Big|\int_{\dist (y, \ell)> 10^4\frac{K_1}K} e^{i\phi(x, y)} f(y) dy\Big| &\leq \sum_{\dist (y_\alpha, \ell)> 10^3 \frac {K_1}K} |T_\alpha
f(x)|< K^2K^{-4} c_*
\\
&<K^{-2} \int|T_{\alpha_*} f(a-z)|\zeta_K(z)dz.
\tag 6.1
\endalign
$$

Hence
$$
\Big|\int_{\dist (y, \ell)> 10^4\frac {K_1}K} e^{i\phi(x, y)} f(y)dy\Big|^q< K^{-2q} \sum^{\sim K^2}_{\alpha=1} \int |T_\alpha f(x-z)|^q \zeta_K(z)dz\tag
6.2
$$
and by (5.5), the corresponding contribution is at most
$$
K^{-2}.K^{\frac 2q}.K^{\frac 4q-2} \ Q_{R/K} < K^{-2} Q_{R/K}.\tag 6.3
$$
Considering the  partition $\Omega=\bigcup\tilde\Omega_\tau$ in balls of size $\frac 1{K_1}$ and fixing $x\in B(a, K)$, there are clearly
the following alternatives

\noindent
(6.4) \ $|Tf(x)|< C \max\limits_\tau\Big|\int_{\tilde\Omega_\tau} e^{i\phi(x, y)} f(y)dy\Big|$.

\noindent
(6.5) \ There are $\tau, \tau'$ such that $\dist (\tilde\Omega_\tau, \tilde\Omega_{\tau'})> \frac {10^6}{K_1}$ and
$$
\Big|\int_{\tilde\Omega_\tau} e^{i\phi(x, y)} f(y)dy\Big|, \Big|\int_{\tilde\Omega_{\tau'}} e^{i\phi(x, y)} f(y) dy\Big|>\frac 1{10K_1^2}|Tf(x)|.
$$

If (6.4), write
$$
|Tf(x)|\leq C\Big[\sum^{\sim K_1^2}_{\tau=1} \Big|\int_{\tilde\Omega_\tau} e^{i\phi(x, y)} f(y)dy\Big|^q\Big]^{\frac 1q}=(6.6)
$$
and by (5.4), (5.5)
$$
\Vert (6.6)\Vert_{L^q(B_R)}\leq  K_1^{\frac 2q}.K_1^{\frac 4q-2} Q_{R/K_1}<  K_1^{-2(1-\frac 3q)} Q_{R/K_1}.\tag 6.7
$$
Assume (6.5). Estimate further
$$
\align
\Big|\int_{\tilde\Omega_\tau} e^{i\phi(x, y)}f(y)dy\Big| & \leq \Big|\sum_{\Sb {\Omega_\alpha\subset\tilde\Omega_\tau} \\ {\dist(y_\alpha, \ell)\leq 10^3\frac
{K_1}K}\endSb} e^{i\phi(x, y_2)} (T_\alpha f)(x)\Big| +\sum_{\Sb {\Omega_\alpha\subset\tilde\Omega_\tau} \\ {\dist (y_\alpha, \ell)> 10^3 \frac {K_1}K}\endSb}
|T_\alpha f|\\
&{}\\
&=(6.8)+(6.9)
\endalign
$$
and similarly for $|\int_{\tilde\Omega_{\tau'}} e^{i\phi(x, y)} f(y) dy|$.

The contribution of (6.9) was evaluated in (6.1), (6.3).

Thus it remains to obtain a bound on
$$
\int_{B(a, K)} \Big|\sum_{\Sb \Omega_\alpha\subset\tilde\Omega_\tau \\ \dist (y_{\alpha}, \ell)\leq 10^{3}\frac{K_1}K\endSb}
e^{i\phi(x, y_\alpha)}(T_\alpha f)(x)\Big|^{\frac q2} \ \Big|\sum_{\Sb \Omega_\alpha \subset \tilde\Omega_{\tau'} \\
{\dist (y_{\alpha}, \ell)\leq 10^3\frac {K_1}K}\endSb} e^{i\phi(x, y_\alpha)}(T_\alpha f) (x)\Big|^{\frac q2} dx \tag 6.10
$$

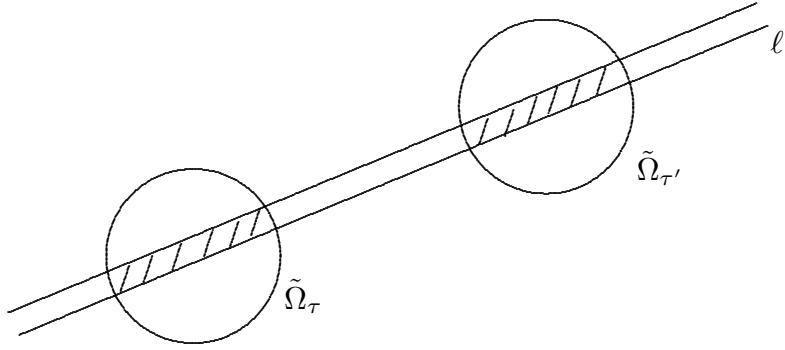
\font\thinlinefont=cmr5
$$
\hbox{\beginpicture
\setcoordinatesystem units <.60000cm,.60000cm>
\setshadesymbol ({\thinlinefont .})
\setlinear
%
% Fig ELLIPSE
%
\linethickness= 0.500pt
\setplotsymbol ({\thinlinefont .})
{\ellipticalarc axes ratio  1.609:1.609  360 degrees 
	from  7.556 15.100 center at  5.948 15.100
}%
%
% Fig ELLIPSE
%
\linethickness= 0.500pt
\setplotsymbol ({\thinlinefont .})
{\ellipticalarc axes ratio  1.609:1.609  360 degrees 
	from 14.069 17.860 center at 12.461 17.860
}%
%
% Fig POLYLINE object
%
\linethickness= 0.500pt
\setplotsymbol ({\thinlinefont .})
{\plot  4.763 14.922  4.604 14.446 /
}%
%
% Fig POLYLINE object
%
\linethickness= 0.500pt
\setplotsymbol ({\thinlinefont .})
{\plot  5.190 15.098  5.031 14.622 /
}%
%
% Fig POLYLINE object
%
\linethickness= 0.500pt
\setplotsymbol ({\thinlinefont .})
{\plot  6.303 15.574  6.145 15.098 /
}%
%
% Fig POLYLINE object
%
\linethickness= 0.500pt
\setplotsymbol ({\thinlinefont .})
{\plot  6.780 15.733  6.621 15.257 /
}%
%
% Fig POLYLINE object
%
\linethickness= 0.500pt
\setplotsymbol ({\thinlinefont .})
{\plot 11.383 17.638 11.225 17.162 /
}%
%
% Fig POLYLINE object
%
\linethickness= 0.500pt
\setplotsymbol ({\thinlinefont .})
{\plot 11.860 17.797 11.701 17.321 /
}%
%
% Fig POLYLINE object
%
\linethickness= 0.500pt
\setplotsymbol ({\thinlinefont .})
{\plot 12.287 18.034 12.129 17.558 /
}%
%
% Fig POLYLINE object
%
\linethickness= 0.500pt
\setplotsymbol ({\thinlinefont .})
{\plot 12.668 18.224 12.510 17.748 /
}%
%
% Fig POLYLINE object
%
\linethickness= 0.500pt
\setplotsymbol ({\thinlinefont .})
{\plot 13.081 18.352 12.922 17.875 /
}%
%
% Fig POLYLINE object
%
\linethickness= 0.500pt
\setplotsymbol ({\thinlinefont .})
{\plot 13.557 18.574 13.399 18.098 /
}%
%
% Fig POLYLINE object
%
\linethickness= 0.500pt
\setplotsymbol ({\thinlinefont .})
{\plot  2.532 14.055 16.343 19.770 /
}%
%
% Fig POLYLINE object
%
\linethickness= 0.500pt
\setplotsymbol ({\thinlinefont .})
{\plot  2.735 13.642 16.546 19.357 /
}%
%
% Fig POLYLINE object
%
\linethickness= 0.500pt
\setplotsymbol ({\thinlinefont .})
{\plot  5.721 15.325  5.563 14.848 /
}%
%
% Fig POLYLINE object
%
\linethickness= 0.500pt
\setplotsymbol ({\thinlinefont .})
{\plot  7.173 15.962  7.015 15.486 /
}%
%
% Fig TEXT object
%
\put{$\ell$} [lB] at 16.599 18.850
%
% Fig TEXT object
%
\put{$\tilde\Omega_{\tau'}$} [lB] at 14.129 16.510
%
% Fig TEXT object
%
\put{$\tilde\Omega_\tau$} [lB] at  7.620 14.129
\linethickness=0pt
\putrectangle corners at  2.506 19.795 and 16.705 13.477
\endpicture}
$$

\vskip .4 true in

By H\"older's inequality, assuming $q< 4$
$$
(6.10) \lesssim K^{3(1-\frac q4)}\Big[\int_{B(a, K)} |\cdots|^2 \ |\cdots|^2 dx\Big]^{q/4}\tag 6.11
$$
Consider
$$
\align
&\int_{B(a, K)} |\cdots|^2 \ |\cdots|^2 \leq\\
&\sum_{\Sb \Omega_{\alpha_1}, \Omega_{\alpha_2} \subset\tilde\Omega_\tau\cap\Delta \\
\Omega_{\alpha_1'}, \Omega_{\alpha_2'} \subset \tilde\Omega_{\tau'} \cap \Delta\endSb}
\Big|\int_{B(a, K)} T_{\alpha_1} f  \ \overline{T_{\alpha_2} f} \ \overline{T_{\alpha_1'} f} \ T_{\alpha_{2}'} f
\quad e^{i[\phi(x, y_{\alpha_1})-
\phi(x,  y_{\alpha_2})\cdots]} dx\Big|
\tag 6.12
\endalign
$$
where $\Delta=\Big\{y\in B(0, 1); \dist (y, \ell)< 10^3 \frac{K_1}K\Big\}$.

Rewriting
$$
\align
&\phi (x, y_{\alpha_1})-\phi(x, y_{\alpha_2}) -\phi(x, y_{\alpha_1'})+\phi(x, y_{\alpha_2'})=\\
&< (x_1, x_2), y_{\alpha_1}-y_{\alpha_2} -y_{\alpha_1'}+ y_{\alpha_2'}>
+ x_3\big(\phi_1(y_{\alpha_1})-\phi_1(y_{\alpha_2})-\phi_1(y_{\alpha_1'})+\phi_1(y_{\alpha_2'})\big)
\endalign
$$
we see that in (6.12) we may restrict the summation to those quadruples $(\alpha_1, \alpha_2, \alpha_1', \alpha_2')$ for which
$$
\qquad\qquad \left\{\aligned
&|y_{\alpha_1}-y_{\alpha_2}- y_{\alpha_1'}+y_{\alpha_2'}|\lesssim \frac 1K \  \qquad\qquad\qquad\qquad\qquad \qquad\qquad\qquad {(6.13)}\\
&|\phi_1(y_{\alpha_1})-\phi_1(y_{\alpha_2})-\phi_1(y_{\alpha_1'})+\phi_1(y_{\alpha_2'})|\lesssim \frac 1K\quad \qquad\qquad\qquad\qquad {(6.13')}
\endaligned
\right.
$$
Let $\ell=b+\mathbb R v$ $(|v|=1)$ and $|y_{\alpha_i}-(b+t_iv)| < 10^3 \frac {K_1}{K}, |y_{\alpha_1'}-(b+ t_i'v)|<10^3\frac{K_1}{K}$.

Recall from (6.5) that
$$
|t_1-t_2|, |t_1'-t_2'|\leq \frac 2{K_1}, |t_1-t_1'|> \frac{10^6}{K_1}.
$$
Hence (6.13), (6.13') imply by the preceding
$$
\qquad \qquad  \left \{\aligned
&|t_1-t_2-t_1'+ t_2'|\lesssim C\frac{K_1}K\qquad\qquad\qquad\qquad\qquad\qquad\qquad \qquad \qquad \ {(6.14)}\\
&|t_1^2-t_2^2 -(t_1')^2 +(t_2')^2|\lesssim C\frac{K_1}K \qquad\qquad\qquad\qquad\qquad\qquad\qquad \quad \ \ {(6.14')}
\endaligned
\right.
$$
and we obtain from the separation property that
$$
|(t_1+t_2)-(t_1'+t_2')|\lesssim C\frac{K_1^2}{K}.\tag 6.14$''$
$$
Hence $|t_1 -t_2|, |t_1'-t_2'|< C\frac{K_1^2}K$, thus $|y_{\alpha_1}- y_{\alpha_2}|, |y_{\alpha_1'}-y_{\alpha_2'}|< C\frac{K_1^2}K$.

Consequently
$$
(6.12) \lesssim K_1^{8} \sum_{\Sb \Omega_\alpha \subset \tilde\Omega_\tau\cap\Delta\\  \Omega_{\alpha'} \subset \tilde\Omega_{\tau'} \cap\Delta\endSb}
\int_{B(a, K)} |(T_\alpha f)(x)|^2 |(T_{\alpha'}f)(x)|^2 dx\tag 6.15
$$
and
$$
\align
(6.10), (6.11) &\lesssim K^{3(1-\frac q4)} K_1^{2q} K^{\frac {3q}4}
 \Big[\sum_{\Omega_\alpha \subset\tilde\Omega_\tau\cap\Delta} c_\alpha^2 \Big]^{\frac q4}\Big[\sum_{\Omega_{\alpha'}\subset
\tilde\Omega_{\tau'} \cap\Delta} c^2_{\alpha'}\Big]^{\frac q4}\\
&\lesssim K^3 K_1^{2q} \Big(\frac K{K_1}\Big)^{(\frac q2-1)} \Big[\sum c_\alpha^{q }\Big]\tag 6.16\\ 
%\qquad\qquad\qquad\qquad\qquad\qquad\qquad\quad\qquad \ {(6.16)}\\
&<K_1^{\frac{3q} 2 +1} K^{\frac q2-1} \sum_\alpha\int\Big[\int_{B(a, K)} |T_\alpha f(x-z)|^q dx\Big] \zeta_K(z) dz. \tag{6.16$'$} 
%\quad\qquad \qquad \quad \ {(6.16')}
\endalign
$$

Summing over the balls $B(a, K)$ implies an estimate
$$
K_1^{\frac 32+\frac 1q} K^{\frac 12-\frac 1q}\Big(\sum_\alpha\Vert T_\alpha f\Vert^q_{L^q(B_R)}\Big)^{\frac 1q} < K_1^{\frac 32+\frac 1q} K^{5/q-3/2}
\ Q_{R/K}.\tag 6.17
$$
Collecting contributions (4.1), (5.7), (6.7), (6.3), (6.17) implies that
$$
Q_R\lesssim C(K)R^\ve+K_1^{-2(1-\frac 3q)} Q_{R/K_1}+ K^{-2} \ Q_{R/K}+K_1^{\frac 32+\frac 1q} K^{\frac 5q-\frac 32}\ Q_{R/K}\tag 6.18
$$
and hence an appropriate choice of $K_1, K$ shows that
$$
Q_R\ll R^\ve \text { for } q > \frac {10}3.\tag 6.19
$$

\noindent
{\bf Remark.} The use of different scales in previous analysis (and even more so in \S3) is reminiscent of the
`induction on scales' approach form [Wo2] and [T1], although the present argument is considerably simpler.
In particular, it suffices to take $K, K_1$ to be large constants, rather than $R$-dependent (i.e.
$R^\ve$-factors), though this point is inessential.

\bigskip

\noindent
{\bf (7).}  One may also consider the hyperbolic case, for instance
$$
\phi(x, y) =x_1y_1+x_2 y_2 +x_3 y_1 y_2.\tag 7.1
$$

The hyperbolic case was studied by Vargas in \cite{V}, adapting the
bilinear method.  She proved the same estimates in the hyperbolic
case that Tao proved in the elliptic case - 
in particular that the restriction operator 
 is bounded from $L^\infty$ into $L^p$ for $p > 10/3$.  Our method
gives nearly the same estimate, losing a factor of $R^\epsilon$.

The preceding may be repeated verbatim, except for the analysis of (6.13').
The condition becomes $(v_1^2+v_2^2=1)$
$$
|v_1| \ |v_2| \ |t_1^2 - t_2^2 -(t_1')^2+(t_1')^2| \lesssim C\frac{K_1}K\tag 7.2
$$
and the case where $v_1$ or $v_2$ is small has to be treated separately.

Suppose $|v_2|<\frac 1{K_1}$.
Let $\Omega=\bigcup_{1\leq s\lesssim K_1} \omega_s$ be a partition in horizontal stripes of width $\frac 1{K_1}$.
Recalling (6.1)-(6.3), for $x\in B(a, R)$, the only significant contribution to $Tf(x)$ is given by
$$
2\max_s\Big|\int_{\omega_s} e^{i\phi(x, y)}f(y) dy\Big|\lesssim \Big[\sum_s \Big|\int_{\omega_s} e^{i\phi(x, y)} f(y) dy\Big|^q\Big]^{\frac 1q}\tag 7.3
$$
since $\ell=b+tv= b+te_1 +0(\frac 1{K_1})$ by assumption on $v$.

The contribution of (7.3) is at most
$$
K_1^{\frac 1q}. \Big\Vert\int_\omega e^{i\phi(x, y)} f(y) dy\Big\Vert_{L^q(B_R)}\tag 7.4
$$
where $\omega=[0, 1]\times [0, \frac 1{K_1}]$.

A rescaling $(x, y)\mapsto (x_1, K_1 x_2, K_1 x_3; {K_1}x_3;  y_1, \frac 1{K_1}y_2)$ shows that
$$
\Big\Vert\int_\omega e^{i\phi(x, y)}f(y) dy\Big\Vert_{L^q(B_R)}\leq K_1^{-1+\frac 2q}Q_R
$$
which in (6.18) gives an extra term $K_1^{-1+\frac 3q} Q_R$.

\bigskip

\noindent
{\bf \S3. Higher Dimensional Restriction Estimates}

The method presented in \S2 easily generalizes to arbitrary dimension, considering the Fourier restriction/extension problem for a smooth, compact hyper-surface $S$ in $\mathbb R^n$ with
positive definite second fundamental form.
For $x\in S$, denote $x' \in S^{(n-1)}$ the normal vector at the point $x$ and let $^\sim: S^{(n-1)} \to S$ be the Gauss map.
Thus $\widetilde {x'} =x$.

In this section, we establish Theorem 1$'$, implying in turn Theorem 1 by the `$\ve$-removal lemma' presented in the Appendix.

\noindent
{\bf 1.} Let $U_1, \ldots, U_n \subset S$ be small caps such that $|x_1'\wedge\cdots \wedge x_n'|>c$ for $x_i\in U_i$.

Let $M$ be large and $\mathcal D_i\subset U_i$ $(1\leq i\leq n)$ discrete sets of $\frac 1M$-separated points.

Let $B_M\subset\mathbb R^n$ be a ball of radius $M$.
Then, for $q=\frac{2n}{n-1}$
$$
\nint_{B_M}\prod^n_{i=1} \Big|\sum_{\xi\in \mathcal D_i} a(\xi) e^{ix.\xi}\Big|^{q/n} \ll M^\ve \prod^n_{i=1} \Big[\sum_{\xi\in\mathcal D_i}
|a(\xi)|^2\Big]^{\frac q{2n}}.
\tag 1.1
$$

\noindent
{\it Proof.}

This is just a discretized version of Theorem 1.16 in [BCT] as our assumption on $U_1, \ldots, U_n$ ensures the required transversality 
condition (see the discussion in the beginning of \S5).

We can assume $B_M$ centered at $0$.
Introduce functions $g_i$ on $U_i$ defined by
$$
\left\{
\aligned
&{g_i(\zeta)= a(\xi) \text { if } |\zeta-\xi|< \frac cM, \xi\in\mathcal D_i}\\
&{g_i(\zeta)= 0 \ \text { otherwise}}
\endaligned\right.
\tag  1.2
$$
($c>0$ a small constant). One may then replace $\sum_{\xi\in\mathcal D_i} a(\xi) e^{ix. \xi}$ by
\hfill\break $c'M^{n-1} \int_S g_i(\zeta) e^{ix.\zeta} \sigma (d\zeta)\text { if } x\in B_M$.
Hence
$$
\aligned
&\int_{B_M} \prod^n_{i=1} \Big|\sum_{\zeta \in\mathcal D_i} a(\xi) e^{ix.\xi}\Big|^{q/n} dx \lesssim\\
&M^{(n-1)q} \int_{B_M} \prod^n_{i=1} \Big|\int_S g_i(\zeta) e^{ix\zeta} \sigma(d\zeta)\Big|^{q/n} dx \overset{\text {[BCT]}} \to \ll\\
&M^{(n-1)q+\ve}\prod^n_{i=1} \Vert g_i\Vert^{q/n}_{L^2(U_i)} \sim M^{\frac {n-1}2 q+\ve} \prod^n_{i=1} \Big[\sum_{\xi\in\mathcal D_i}|a(\xi)|^2\Big]^{\frac
q{2n}}.
\endaligned\tag1.3
$$
Since $\nint_{B_M}$ refers to the average, (1.1) follows, since $q=\frac {2n}{n-1}$.
\bigskip

\noindent
{\bf 2.} Let $S\subset\mathbb R^n$ be as above and $2\leq m\leq n$.
Let $V$ be an $m$-dimensional subspace of $\mathbb R^n$, $P_1, \ldots, P_m\in S$ such that
$$
P_1', \ldots, P_m' \in V \text { and } \ |P_1'\wedge\cdots\wedge P_m'|>c \tag 2.1
$$
and $U_1, \ldots, U_m \subset S$ sufficiently small neighborhoods of $P_1, \ldots, P_m$.

Let $M$ be large and $\mathcal D_i\subset U_i$ $(1\leq i\leq m)$ discrete sets of $\frac 1M$-separated points 
$\xi\in S$ such that $\dist(\xi', V)<\frac cM$.
Let $g_i\in L^\infty (U_i)(1\leq i\leq m)$.
Then letting $q=\frac {2m}{m-1}$
$$
\aligned
&\nint_{B_M} \prod^m_{i=1} \Big|\sum_{\xi\in\mathcal D_i}\Big(\int_{|\zeta -\xi|<\frac cM} g_i(\zeta) e^{ix.\zeta}
\sigma(d\zeta)\Big)\Big|^{q/m} dx\ll\\
&M^\ve\Big\{\nint_{B_M} \prod^m_{i=1} \Big[\sum_{\xi\in\mathcal D_i} \Big|\int_{|\zeta-\xi|<\frac cM} g_i(\zeta) e^{ix.\zeta}
\sigma(d\zeta)\Big|^2\Big]^{1/2m}\Big\}^q.
\endaligned
\tag 2.2
$$
\bigskip

\noindent
{\it Proof.}

Performing a rotation, we may assume $V=[e_1, \ldots, e_m]$ and denote $\tilde V$ the image of $V\cap S^{(n-1)}$ under the Gauss map.
Let again $B_M$ be centered at $0$.
For each $\xi\in \bigcup^m_{i=1}\mathcal D_i$ there is by assumption some $\hat\xi\in  S\cap \tilde V, |\xi-\hat\xi|< \frac cM$.
Write
$$
\int_{|\zeta-\xi|<\frac cM} g_i(\zeta) e^{ix.\zeta} \sigma(d\zeta)=e^{i x. \hat\xi} \int_{|\zeta-\xi|<\frac cM}
g_i(\zeta) e^{ix.(\zeta-\hat\xi)} \sigma (d\zeta).\tag 2.3
$$
Since in the second factor of (2.3), $|\zeta-\hat\xi|=o(\frac 1M)$, we may view it as a constant $a(\xi)$ on $B_M\subset \mathbb R^n$.

Thus we need to estimate
$$
\nint_{B_M} \Big\{\prod^m_{i=1} \Big|\sum_{\xi \in \Cal D_i} e^{ix.\hat\xi} a(\xi)\Big|^{q/m}\Big\} dx.\tag 2.4
$$
Writing $x=(u, v) \in B_M^{(m)} \times B_M^{(n-m)}$, (2.4) may be bounded by
$$
\max_{v\in B_M^{(n-m)}} \nint_{B_M^{(m)}} \Big\{\prod^m_{i=1} \Big|\sum_{\xi\in\Cal D_i} e^{iu.\pi_m(\hat\xi)} a_v(\xi)\Big|^{q/m}\Big\} du\tag 2.5
$$
with $a_v(\xi)=e^{iv.\hat\xi} a(\xi)$.

Since $S$ has positive definite second fundamental form, $\pi_m(S\cap \tilde V)\subset V= [e_1, \ldots, e_m]$ is a hypersurface in $V$ with 
same property and the normal vector at $\pi_m(\hat\xi)= (\hat\xi)' \in V$.
Since (2.1), application of (1.1) with $n$ replaced by $m$ and $\Cal D_i$ by $\{\pi_m\hat\xi; \xi \in\Cal D_i\}$ gives the estimate on (2.5)
$$
\ll M^\ve \prod^m_{i=1} \Big[\sum_{\xi\in \Cal D_i} |a(\xi)|^2\Big]^{q/2m}
$$
and (2.2) follows.

\noindent
{\bf 3.} Essential use is made of scaling.

Denote $Q_R^{(p)}$ a bound on
$$
\Big\Vert \int_S g(\xi) e^{ix.\xi} \sigma(d\xi)\Big\Vert_{L^p (B_R)}
$$
with $g\in L^\infty (S), |g|\leq 1$ and with $S$ as specified in the beginning of \S3.

Parametrize $S$ (locally) as
$$
\left\{
\aligned
&\xi_i=y_i \quad (1\leq i\leq n-1)\\
&\xi_n= y_1^2+\cdots+ y^2_{n-1} + O(|y|^3)\endaligned
\right.\tag 3.1
$$
with $y$ taken in a small neighborhood of $0$.

Let $U_\rho$ be a $\rho$-cap on $S$ and evaluate
$$
\Big\Vert\int_{U_\rho} g(\xi) e^{ix.\xi} \sigma(d\xi)\Big\Vert_{L^p(B_R)}.
$$
Thus in (3.1) we restrict $y$ to a ball $B(a, \rho)\subset\mathbb R^{n-1}$ and evaluate
$$
\Big\Vert\int_{B(a, \rho)} g(y) e^{i[x_1 y_1+\cdots+ x_{n-1} y_{n-1}+x_n(|y|^2+O(|y|^3))]} dy\Big\Vert_{L^p(B_R)}.\tag 3.2
$$
A shift $y\mapsto y-a$ and a change of variables $x_i'=x_i+x_n(2a_i+\cdots)(1\leq i< n)$ permits to set $a=0$.
Rescale $y=\rho y'$ to obtain
$$
\rho^{n-1} \Big\Vert\int_{B(0, 1)} g(\rho y') e^{i[\rho x_1y_1'+\cdots+\rho x_{n-1} y_{n-1}' 
+\rho^2x_n(|y'|^2+\rho O(|y'|^{3}))]} dy'
\Big\Vert_{L^p(B_R)}
$$
and a further rescaling in $x$, $x_i' =\rho x_i (1\leq i\leq n-1),  x_n'=\rho^2 x_n$, gives
$$
\aligned
&\rho^{n-1-(n+1)/p}\Big\Vert \int_{B(0, 1)} g(\rho y') e^{i[x_1' y_1'+\cdots+ x_{n-1}' y_{n-1}'+x_n'(|y'|^2+\rho 
O(|y'|^3))]}
dy'\Big\Vert_{L^p(B_{\rho R})}\\
& \leq \rho^{n-1-(n+1)/p}  \ Q_{\rho R}^{(p)}
\endaligned
\tag 3.3
$$

\noindent
{\bf 4.} Let $g\in L^\infty(S), |g|\leq 1$ and consider for $x\in B_R$
$$
\int_Sg(\xi) e^{ix.\xi} \sigma(d\xi).\tag 4.1
$$
Let
$$
R^\ve\gg K_n\gg K_{n-1}\gg\cdots\gg K_1
$$
be suitably chosen.

Start decomposing $S=\bigcup_\alpha U_\alpha (\frac 1{K_n})$ in caps of size $\frac 1{K_n}$ and write
$$
(4.1) =\sum_\alpha\int_{U_\alpha(\frac 1{K_n})} g(\xi) e^{ix.\xi}\sigma (d\xi)=\sum_\alpha c_\alpha(x).
$$
Fixing $x$, there are 2 possibilities

(4.2) There are $\alpha_1, \alpha_2, \ldots, \alpha_n$ such that
$$
|c_{\alpha_1}(x)|, \ldots, |c_{\alpha_n}(x)|> K_n^{-n}\max_\alpha |c_\alpha(x)|\tag 4.3
$$
and
$$
|\xi_1'\wedge\cdots\wedge \xi_n'|> c(K_n) \text { for } \xi_i\in U_{\alpha_i}.
\tag 4.4
$$

(4.5) The negation of (4.2), which implies that there is an $(n-1)$-dim subspace $V_{n-1}$ such that
$$
|c_\alpha(x)|\leq K_n^{-n} \max_\alpha |c_\alpha(x)|\ \text { if } \ \dist (U_\alpha, \tilde V_{n-1})\gtrsim \frac 1{K_n}.
$$
If (4.2), clearly by (4.3)
$$
\Big|\int_Sg(\xi) e^{ix.\xi} \sigma (d\xi)\Big| \leq K_n^{n-1} \max |c_\alpha(x)| \leq K_n^{2n-1} \Big[\prod^n_{i=1} |c_{\alpha_i}(x)|\Big]^{\frac 1n}
$$
and
$$
\int_{x (4.2)} \Big|\int_Sg(\xi) e^{ix.\xi} \sigma(d\xi)\Big|^p \lesssim K_n^{p(2n-1)} \sum_{\Sb {\alpha_1, \ldots, \alpha_n}\\{(4.4)}\endSb}
\int_{B_R}\prod^n_{i=1} \Big|\int_{U_{\alpha_i(\frac 1{K_n})}} g(\xi) e^{ix.\xi}\sigma(d\xi)\Big|^{\frac pn}.
\tag 4.6
$$
In view of (4.4), the [BCT]-estimate applies to each (4.6) term.
Thus
$$
\int_{B_R} \prod^n_{i=1} \Big|\int_{U_{\alpha_i(\frac 1{K_n})}} g(\xi) e^{ix.\xi} \sigma (d\xi)\Big|^{\frac 2{n-1}} dx\ll C(K_n)R^\ve.
\tag 4.7
$$
Assuming
$$
p\geq \frac{2n}{n-1}\tag 4.8
$$
we see that
$$
(4.6) < C(K_n)R^\ve\tag 4.9
$$
(here and in the sequel, $C(K)$ refers to some power of $K$.)

Next consider the case (4.5). Thus
$$
\aligned
|(4.1)|&\leq \Big|\int_{\dist(\xi, \tilde V_{n-1})\lesssim \frac 1{K_n}} g(\xi) e^{ix.\xi}\sigma(d\xi)\Big|+\frac 1{K_n} \max_\alpha\Big|\int_{U_\alpha(\frac
1{K_n})} g(\xi) e^{ix.\xi}\sigma(d\xi)\Big|\\
& = (4.10)+(4.11)
\endaligned
$$
where $V_{n-1}$ depends on $x$.

Note that, using the argument explained earlier in \S1, we may view $|c_\alpha(x)|$ as essentially constant on balls of size $K_n$
(literally speaking, this is of course incorrect and what was done is a replacement of $|c_\alpha(x)|$ by a majorant $|c_\alpha|* \eta_{_{K_n}},
\eta_{_K}(x)=\frac 1{K^d}\eta(\frac xK)$ and $\eta$ a suitable bump-function -- we do not repeat these technicalities here.)

Thus the bound (4.10) + (4.11) may be considered valid on $B(x, K_n)$, with a same linear space $V_{n-1}$.

The contribution of (4.11) to $\Vert\int g(\xi) e^{ix.\xi} \sigma (d\xi)\Vert_p$ is bounded by
$$
\aligned
\frac 1{K_n} \Big(\sum_\alpha \Big\Vert\int_{U_\alpha} g(\xi) e^{ix. \xi} \sigma(d\xi) \Big\Vert^p_p\Big)^{\frac 1p} &\lesssim\frac 1{K_n}.
K_n^{\frac{n-1}p}.\Big(\frac 1{K_n}\Big)^{n-1-\frac {n+1} p} Q^{(p)}_{R/K_n}\\
&=\Big(\frac 1{K_n}\Big)^{n(1-\frac 2p)} Q^{(p)}_{{R/K_n}}<\frac 1{K_n} Q_R^{(p)}.
\endaligned
\tag4.12
$$

Consider the term (4.10).
Proceeding similarly, write for $x\in B(\bar x, K_n)$
$$
\aligned
\int_{\dist (\xi, \tilde V_{n-1})\lesssim \frac 1{K_n}} g(\xi) e^{ix.\xi} \sigma (d\xi)&=\\
\sum_\alpha\int_{U _\alpha (\frac 1{K_{n-1}})\cap[\dist|\xi, \tilde V_{n-1})\lesssim \frac 1{K_n}]} &g(\xi) e^{ix.\xi}\sigma(d\xi) =\sum_\alpha
c_\alpha^{(n-1)}(x).
\endaligned
\tag 4.13
$$
We distinguish the cases

(4.14) There are $\alpha_1, \ldots, \alpha_{n-1}$ such that
$$
|c_{\alpha_1}^{(n-1)} (x)|, \ldots, |c^{(n-1)}_{\alpha_{n-1}} (x)|> K_{n-1}^{-(n-1)}\max_\alpha|c_\alpha^{(n-1)}(x)|\tag 4.15
$$
and
$$
|\xi_1'\wedge\ldots\wedge\xi_{n-1}'|> c(K_{n-1}) \ \text { for } \ \xi_i\in U_{\alpha_i} \Big(\frac 1{K_{n-1}}\Big).\tag 4.16
$$

(4.17) Negation of (4.14), implying that there is an $(n-2)$-dim subspace $V_{n-2}\subset V_{n-1}$ (depending on $x$) such that
$$
|c_\alpha^{(n-1)}(x)|< K_{n-1}^{-(n-1)} \max_\alpha |c_\alpha^{(n-1)} (x)| \text { for } \dist (U_\alpha, \tilde V_{n-2})\gtrsim \frac 1{K_{n-1}}.
$$
This space $V_{n-2}$ can then again be taken the same on a $ K_{n-1}$-neighborhood of $x$.

We analyze the contribution of (4.14).  By (4.15)
$$
|(4.13)| < K_{n-1}^{2n-3} \Big[\prod^{n-1}_{i=1} |c_{\alpha_i}^{(n-1)} (x) |\Big]^{\frac 1{n-1}}\tag 4.18
$$
and hence
$$
\aligned
&\operatornamewithlimits\nint\limits_{\Sb{B(\bar x, K_n)}\\ {x\text { satisfies } (4.14)}\endSb}\Big|\int\limits_{\dist(\xi, \tilde
V_{n-1})\lesssim \frac 1{K_n}} g(\xi) e^{ix.\xi}
\sigma(d\xi)\Big|^p\leq\\
&K_{n-1}^{p(2n-3)}\sum_{\Sb {\alpha_1, \ldots, \alpha_{n-1}}\\ {(4.16)}\endSb} \nint_{B(\bar x, K_n)}\Big\{\prod^{n-1}_{i=1} \Big|\int_{U_{\alpha_i}
(\frac 1{K_{n-1}}) \cap [\dist(\xi, \tilde V_{n-1})\lesssim \frac 1{K_n}]} g(\xi) e^{ix.\xi} \sigma(d\xi)\Big|^{p/n-1}\Big\}
\endaligned
\tag 4.19
$$
We use the bound (2.2) to estimate the individual integrals
$$
(4.20) \ \nint_{B(\bar x, K_n)} \Big\{\prod^{n-1}_{i=1} \Big|\int_{U_{\alpha_i} (\frac 1{K_{n-1}}) \cap [\dist(\xi, \tilde V_{n-1})\lesssim \frac 1{K_n}]}
g(\xi)e^{ix.\xi} \sigma(d\xi) \Big|\Big\}^{\frac q{n-1}} \text { with } q=\frac{2(n-1)}{n-2}.
$$
Thus $m=n-1, V=V_{n-1}$ and $P_i$ is the center of $U_{\alpha_i}(\frac 1{K_{n-1}})$.
Let $M=K_n$ and $\mathcal D_i$ the centers of a cover of $U_{\alpha_i}(\frac 1{K_{n-1}})$ by caps $U_\alpha(\frac1{K_n})$.

By (2.2) we get an estimate
$$
(4.20) \ll K_n^\ve C(K_{n-1}) \Big\{\nint_{B(\bar x, K_n)} \prod^{n-1}_{i=1} \Big[\sum_\alpha^{(i)} \Big|\int_{U_\alpha(\frac 1{K_n})}
g(\xi) e^{ix.\xi} \sigma(d\xi)\Big|^2\Big]^{\frac 1{2(n-1)}}\Big\}^q\tag 4.21
$$
where in $\sum^{(i)}$ the sum is over those $\alpha$ such that $U_\alpha(\frac 1{K_n})\subset U_{\alpha_i}(\frac 1{K_{n-1}})$ and
\hfill\break $U_\alpha(\frac 1{K_n})\cap \tilde V_{n-1} \not=\phi$.
Clearly
$$
(4.21)\ll K_n^\ve \ C(K_{n-1})\Big\{\nint_{B(\bar x, K_n)} \Big[\sum_{U_\alpha(\frac 1{K_n})\cap \tilde V_{n-1} \not=\phi}\Big|\int_{U_\alpha(\frac 1{K_n})}
g(\xi) e^{ix.\xi}\sigma(d\xi)\Big|^2 \Big]^{\frac 12}\Big\}^q.\tag 4.22
$$
If
$$
p\geq \frac{2(n-1)}{n-2} =q \tag 4.23
$$
the contribution of (4.15) may be estimated replacing $p$ by $q=\frac {2(n-1)}{n-2}$, and using the [BCT] bound (4.7) with $n$ replaced by $n-1$ and $K_n$ by
$K_{n-1}$.
This gives a bound $R^\ve$.

Thus we assume
$$
p<\frac{2(n-1)}{n-2}.\tag 4.24
$$
Then
$$
(4.19)^{1/p}\ll C(K_{n-1}) K_n^\ve\nint_{B(\bar x, K_n)} \Big[\sum_{U_a (\frac 1{K_n}) \cap \tilde V_{n-1}\not=\phi} \Big|\int_{U_\alpha(\frac 1{K_n})}
g(\xi) e^{ix. \xi} \sigma(d\xi)\Big|^2\Big]^{\frac 12}.\tag 4.25
$$
Note that $U_\alpha(\frac 1{K_n})\cap \tilde V_{n-1} \not=\phi$ for $\sim K_n^{n-2}$ values of $\alpha$.

Hence, by H\"older's inequality, the integrand in (4.25) is at most
$$
K_n^{(n-2)(\frac 12-\frac 1p)} \Big[\sum_\alpha \Big|\int_{U_\alpha(\frac 1{K_n})} g(\xi)e^{ix.\xi}\sigma(d\xi)\big|^p\Big]^{\frac 1p}\tag 4.26
$$
where $\alpha$ is unrestricted in the $\alpha$-summation.
Substituting (4.26) in (4.25) gives
$$
(4.19)\ll C(K_{n-1})K_n^{(n-2)(\frac p2-1)+\ve}\nint_{B(\bar x, K_n)} \Big[\sum_\alpha\Big|\int_{U_\alpha(\frac 1{K_n})} g(\xi) e^{ix.\xi}
\sigma(d\xi)\Big|^p\Big]
$$
and integrating over $B_R$ permits to bound the (4.14)-contribution by
$$
C(K_{n-1})K_n^{(n-2)(\frac 12-\frac 1p)+\ve}\Big[\sum_\alpha\Big\Vert\int_{U_\alpha(\frac 1{K_n})} g(\xi) e^{ix.\xi} \sigma(d\xi)
\Big\Vert^p_{L^p(B_R)}\Big]^{1/p}.\tag 4.27
$$

Invoking again the rescaling inequality (3.3), this gives
$$
C(K_{n-1}) K_n^{(n-2)(\frac 12-\frac 1p)+\frac{n-1}p - (n-1)+\frac{n+1}p +\ve } Q_{R/K_n} =C(K_{n-1}) K_n^{\frac{n+2}p -\frac n2+\ve}.
\tag 4.28
$$

Taking $K_n$ sufficiently large compared with $K_{n-1}$, we see that the (4.14)-contribution is taken care of if either $p\geq \frac{2(n-1)}{n-2}$ or
$$
p>2+\frac 4n.\tag 4.29
$$
Thus we impose
$$
p>\min\Big(\frac{2(n-1)}{n-2}, \frac{2(n+2)}n\Big).\tag 4.30
$$
Next we need to consider the contribution of (4.17).

The analysis is analogous to the preceding, replacing $n-1$ by $n-2$ and $K_n$ by$K_{n-1}$.
More precisely, if
$$
p<\frac{2(n-2)}{n-3}\tag 4.31
$$
the local estimate (4.25) becomes
$$
c(K_{n-2})K^\ve_{n-1} \nint_{B(\bar x, K_{n-1})} \Big[\sum_{U_\alpha(\frac 1{K_{n-1}})\cap \tilde V_{n-2}\not=\phi}\Big|\int_{U_\alpha(\frac 1{K_{n-1})}}
g(\xi) e^{ix.\xi}\sigma(d\xi)\Big|^2 \Big]^{\frac 12}\tag 4.32
$$
and $U_\alpha(\frac 1{K_{n-1}})\cap \tilde V_{n-2}\not=\phi$ for $\sim K_{n-1}^{n-3}$ values of $\alpha$.

This leads to the condition on $p$
$$
p>\min\Big(\frac{2(n-2)}{n-3}, \frac{2(n+3)}{n+1}\Big).\tag 4.33
$$

The continuation of the process is clear.

\medskip

Eventually we see that the exponent $p$ needs to satisfy
$$
p>2\min \Big\{\frac k{k-1}, \frac {2n-k+1}{2n-k-1}\Big\} \text { for all } 2\leq k\leq n.\tag 4.34
$$
Hence we obtain.

\proclaim
{\bf Theorem 1'}

$Q_R^{(p)} \ll R^\ve$ provided
$$
p\geq 2\frac{4n+3}{4n-3}\quad \text { if } n\equiv 0 \ (\mod 3)
$$
$$
p\geq \frac{2n+1}{n-1}\quad \, \text { if } n\equiv 1 \ (\mod 3)
$$
$$
p\geq \frac {4(n+1)}{2n-1} \quad \text { if } n\equiv 2 \ (\mod 3).
$$
\endproclaim

\bigskip

\noindent
{\bf \S4 Improving Upon the Exponent in the 3D Restriction Problem}

We consider the case of the paraboloid (though the argument generalizes).

Going back to the analysis in \S2, the main idea is to collect the contributions obtained
at different scales, rather than performing an induction on scale argument.
This will allow us to bring into play also T.~Wolff's $\frac 52$-bound for the Kakeya maximal function.
(see [Wo1]).

\bigskip

\noindent
{\bf 1. Representation at scale 1}

\
Fix large parameters $K\gg K_1\gg 1$

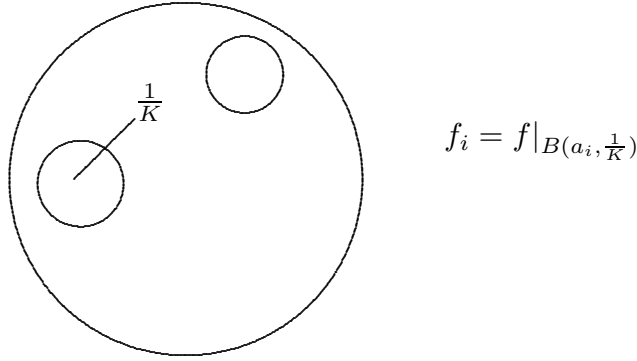
\font\thinlinefont=cmr5
$$
\hbox{\beginpicture
\setcoordinatesystem units <.60000cm,.60000cm>
%
% Fig ELLIPSE
%
\linethickness= 0.500pt
\setplotsymbol ({\thinlinefont .})
{\ellipticalarc axes ratio  0.794:0.794  360 degrees 
	from  6.312 21.188 center at  5.518 21.188
}%
%
% Fig ELLIPSE
%
\linethickness= 0.500pt
\setplotsymbol ({\thinlinefont .})
{\ellipticalarc axes ratio  0.709:0.709  360 degrees 
	from  9.256 23.203 center at  8.547 23.203
}%
%
% Fig ELLIPSE
%
\linethickness= 0.500pt
\setplotsymbol ({\thinlinefont .})
{\ellipticalarc axes ratio  3.258:3.258  360 degrees 
	from 10.719 21.273 center at  7.461 21.273
}%
%
% Fig POLYLINE object
%
\linethickness= 0.500pt
\setplotsymbol ({\thinlinefont .})
{\plot  5.397 21.273  6.509 22.384 /
}%
%
% Fig TEXT object
%
\put{$\frac 1K$ } [lB] at  6.509 22.543
%
% Fig TEXT object
%
\put{$f_i= f|_{B(a_i, \frac 1K)}$} [lB] at 12.224 21.907
\linethickness=0pt
\putrectangle corners at  4.187 24.547 and 15.348 17.998
\endpicture}
$$

Recalling the analysis in \S2, we have
$$
\aligned
|Tf|&\leq C(K) \max_{\Sb {i_1, i_2, i_3}\\ {\text \rm non-collinear}\endSb} (|Tf_{i_1}|.|Tf_{i_2}| \ | Tf_{i_3}|)^{\frac 13} +\max_{\Sb{\mathcal L}\\
{\dist (\mathcal L', \mathcal L'')>\frac 1{K_1}}\endSb} \Big|\sum_{i\in\mathcal L'} Tf_i\Big|^{\frac 12}\Big|\sum_{i\in\mathcal L''} Tf_i\Big|^{\frac 12}
\\
& +\max_a\big| T(f|_{B(a, \frac 1{K_1})})\big|\\
{}\\
&= (1.1)+ (1.2)+(1.3).\\
\endaligned
$$
Here $\mathcal L', \mathcal L'' \subset\mathcal L$ are separated segments of a `line' $\mathcal L$.
\font\thinlinefont=cmr5
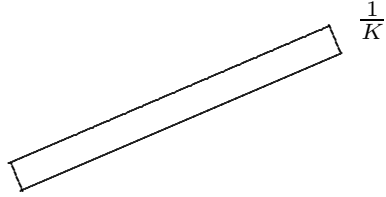
$$
\hbox{\beginpicture
\setcoordinatesystem units <.80000cm,.80000cm>
%
% Fig POLYLINE object
%
\linethickness= 0.500pt
\setplotsymbol ({\thinlinefont .})
{\plot  1.905 22.225  6.350 24.130 /
}%
%
% Fig POLYLINE object
%
\linethickness= 0.500pt
\setplotsymbol ({\thinlinefont .})
{\plot  2.070 21.846  6.515 23.751 /
}%
%
% Fig POLYLINE object
%
\linethickness= 0.500pt
\setplotsymbol ({\thinlinefont .})
{\plot  6.342 24.145  6.500 23.764 /
}%
%
% Fig POLYLINE object
%
\linethickness= 0.500pt
\setplotsymbol ({\thinlinefont .})
{\plot  1.939 22.225  2.098 21.844 /
}%
%
% Fig TEXT object
%
\put{$\frac 1K$} [lB] at  6.731 24.098
\linethickness=0pt
\putrectangle corners at  1.880 24.384 and  7.120 21.819
\endpicture}
$$

Since
$$
\Big[\nint_{B(a, K)} (1.2)^4 \Big]^{\frac 14}\leq C(K_1)\Big(\sum_{i\in\mathcal L} |Tf_i|^2\Big)^{\frac 12}
$$
we may write
$$
(1.2) =\phi.\Big(\sum_{i\in\mathcal L}|Tf_i|^2\Big)^{\frac 12}
$$
with
$$
\Big(\nint_{B(a, K)}|\phi|^4\Big)^{\frac 14} < c(K_1)
$$
and $\phi$ constant on balls of radius 1.

In what follows, we identify small discs $\subset\Omega$ and the corresponding caps $\subset S$ obtained as image under the map
$y\mapsto (y_1, y_2, y_1^2+ y_2^2)$, which are both denoted by $\tau$.

\bigskip

\noindent
{\bf 2. Representation of $Tf_\tau$} (by rescaling).

Let $\tau$ be a $\delta$-cap and rescale.

Up to linear transformation of the form
$$
\left\{
\aligned
x_1' &= x_1+ a_1x_3\\
x_2'&= x_2+ a_2x_3\\
x_3' &=x_3
\endaligned
\right.
$$
and reduction to scale 1 by transformation
$$
\left\{
\aligned
x_1'&= \delta x_1\\
x_2'&= \delta x_2\\
x_3'&=\delta^2 x_3
\endaligned
\right.
$$
we obtain

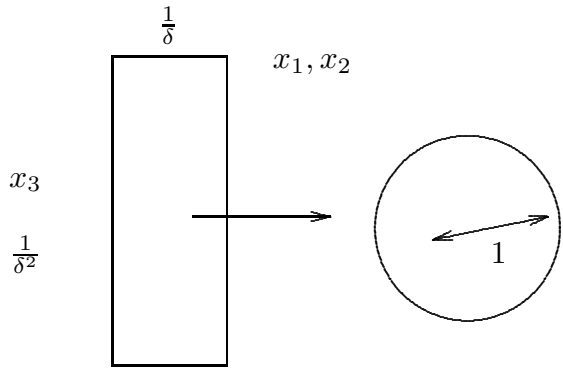
\font\thinlinefont=cmr5
$$
\hbox{\beginpicture
\setcoordinatesystem units <.80000cm,.80000cm>
%
% Fig ELLIPSE
%
\linethickness= 0.500pt
\setplotsymbol ({\thinlinefont .})
{\ellipticalarc axes ratio  1.281:1.281  360 degrees 
	from  8.107 22.225 center at  6.826 22.225
}%
%
% Fig POLYLINE object
%
\linethickness= 0.500pt
\setplotsymbol ({\thinlinefont .})
{\putrectangle corners at  1.905 24.606 and  3.493 20.320
}%
%
% Fig POLYLINE object
%
\linethickness= 0.500pt
\setplotsymbol ({\thinlinefont .})
{\putrule from  3.016 22.384 to  4.921 22.384
%
% arrow head
%
\plot  4.667 22.320  4.921 22.384  4.667 22.447 /
}%
%
% Fig POLYLINE object
%
\linethickness= 0.500pt
\setplotsymbol ({\thinlinefont .})
% arrow head
%
\plot  6.587 22.178  6.350 22.066  6.612 22.054 /
\plot  6.350 22.066  7.938 22.384 /
%
% arrow head
%
\plot  7.701 22.272  7.938 22.384  7.676 22.396 /
%
% Fig TEXT object
%
\put{$\frac 1 \delta$} [lB] at  2.540 24.924
%
% Fig TEXT object
%
\put{$x_1, x_2$} [lB] at  4.128 24.448
%
% Fig TEXT object
%
\put{$x_3$} [lB] at 0.476 22.800

\put{$\frac 1{\delta^2}$} [lB] at  0.476 21.749
%
% Fig TEXT object
%
\put{$1$} [lB] at  7.144 21.749
\linethickness=0pt
\putrectangle corners at  0.476 25.210 and  8.124 20.295
\endpicture}
$$

Applying at unit scale the representation from (1) and scaling back, we obtain
on the $(\frac K\delta \times \frac K\delta\times\frac K{\delta^2})$-box
\bigskip
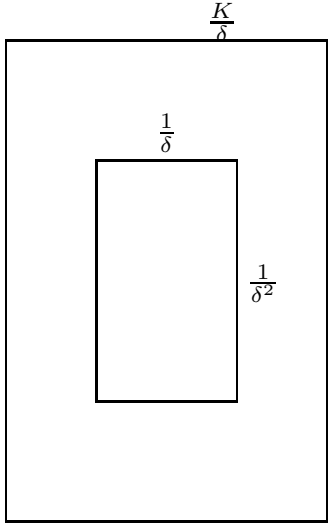
\font\thinlinefont=cmr5
$$
\hbox{\beginpicture
\setcoordinatesystem units <.70000cm,.70000cm>
%
% Fig POLYLINE object
%
\linethickness= 0.500pt
\setplotsymbol ({\thinlinefont .})
{\putrectangle corners at  1.270 25.400 and  6.350 17.780 }%
%
% Fig POLYLINE object
%
\linethickness= 0.500pt
\setplotsymbol ({\thinlinefont .})
{\putrectangle corners at  2.699 23.495 and  4.921 19.685 }%
%
% Fig TEXT object
%
\put{$\frac K\delta$} [lB] at  4.445 25.559
%
% Fig TEXT object
%
\put{$\frac 1\delta$} [lB] at  3.651 23.812
%
% Fig TEXT object
%
\put{$\frac K{\delta^2}$} [lB] at  6.668 21.431
%
% Fig TEXT object
%
\put{$\frac 1{\delta^2}$} [lB] at  5.080 21.431
%
% Fig TEXT object
%
\put{ $_{|Tf_\tau|\leq C(K)  \max\limits_{\Sb{\tau_1, \tau_2, \tau_3}\\ {\text{non-collinear}}\endSb}  |Tf{\tau_1}|^{\frac 13} |Tf\tau_2|^{\frac 13}
.|Tf\tau_3|^{\frac 13}}\qquad  \,  \quad (2.1)$} [lB] at  7.699 24.765
%
% Fig TEXT object
%
\put{$ _{+\phi_\tau \max\limits_{\Cal L} \Big(\sum\limits_{\tau_i\in \Cal L} |Tf\tau_i|^2\Big)^{\frac 12}
 \text { where  $\tau_i$ is a } \frac \delta K-\text{cap}}
\qquad  \   (2.2)$} [lB] at  8.700 22.447
%
% Fig TEXT object
\put{$_{+ \max\limits_{\Sb{\tau' \subset \tau}\\ {\delta/K_1-\text{cap}}\endSb} |Tf\tau'|}\qquad\qquad\qquad\qquad\qquad\qquad
\ \ {(2.3)}$} [lB] at 8.700 19.971
\linethickness=0pt
\putrectangle corners at  1.245 25.845 and 16.569 17.755
\endpicture}
$$

Given a $\delta$-cap $\tau$, denote $\overset o \to\tau$ the polar set

\bigskip
\font\thinlinefont=cmr5
$$
\hbox{\beginpicture
\setcoordinatesystem units <.60000cm,.60000cm>
%
% Fig POLYLINE object
%
\linethickness= 0.500pt
\setplotsymbol ({\thinlinefont .})
{\plot 13.335 19.685 19.050 25.400 /
}%
%
% Fig POLYLINE object
%
\linethickness= 0.500pt
\setplotsymbol ({\thinlinefont .})
{\plot 19.050 25.400 20.955 23.971 /
}%
%
% Fig POLYLINE object
%
\linethickness= 0.500pt
\setplotsymbol ({\thinlinefont .})
{\plot 15.287 18.241 21.002 23.956 / }
% % % Fig POLYLINE object % \linethickness= 0.500pt 
\setplotsymbol ({\thinlinefont .}) 
{\plot 13.339 19.689 15.244 18.260 / }% % % 
%Fig POLYLINE object % \linethickness= 0.500pt 
\setplotsymbol ({\thinlinefont .}) 
{\plot 16.986 22.511 16.510 22.035 / }% % 
% Fig POLYLINE object % \linethickness= 0.500pt 
\setplotsymbol ({\thinlinefont .}) 
{\plot 16.986 22.479 17.939 21.527 / }% % 
% Fig POLYLINE object % \linethickness= 0.500pt 
\setplotsymbol ({\thinlinefont .}) 
{\plot 16.525 22.020 17.477 21.067 / }% % 
% Fig POLYLINE object % \linethickness= 0.500pt 
\setplotsymbol ({\thinlinefont .}) 
{\plot 17.939 21.590 17.462 21.114 / }% %
% Fig TEXT object
%
\put{$\tau\to \overset o \to \tau = (\frac 1\delta\times \frac 1\delta\times\frac 1{\delta^2})\text{ box}$}
 [lB] at  4.477 23.463
%
% Fig TEXT object
%
\put{$\delta$} [lB] at 17.558 22.225
%
% Fig TEXT object
%
\put{$\frac 1\delta$} [lB] at 20.288 24.892
\put{$\overset o\to \tau$} [lB] at 13.652 18.092
%
% Fig TEXT object
%
\put{$\tau$} [lB] at 16.669 21.273
%
% Fig TEXT object
%
\put{$\delta^2$} [lB] at 17.721 20.900
%
% Fig TEXT object
%
\put{$\frac 1{\delta^2}$}  [lB] at 15.558 22.760
\linethickness=0pt
\putrectangle corners at  4.477 25.425 and 21.027 18.002
\endpicture}
$$

\bigskip

On every $K\overset o\to\tau $-box $B$, $\phi_\tau$ satisfies
$$
\aligned
\nint_B\phi_\tau^4 &= \frac 1{|B|} \int_B\phi_\tau^4\\
&= \frac {\delta^4}{K^3}. \delta^{-4}\int_{B(a, K)} \phi_\tau(\delta^{-1} x_1', \delta^{-1}{x_2'}, \delta^{-2} x_3')^4 \, dx_1' dx_2' dx_3'\\
&< C(K_1)
\endaligned
\tag2.4
$$
and $\phi_\tau$ is essentially  constant on $\overset o\to\tau$-boxes.

\bigskip
\noindent
{\bf 3. Iteration}

Apply the decomposition (2.1)-(2.3) to each $Tf_{\tau_i}$ in (2.2) and $Tf_{\tau'}$ in (2.3).

Considering $Tf_{\tau_i}$, let $\phi_{\tau_i}$ be the corresponding factor appearing in (2.2).

Thus $\phi_{\tau_i}$ is constant on $\overset o\to\tau_i$-boxes and $\nint_{B'} \phi^4_{\tau_i} < C(K_1)$ if $B'$ is a $K\overset o\to\tau_i$-box.

Let $B'$ be a $K\overset o\to\tau_i$-box and subdivide $B'$ as
$$
B'=\bigcup B_\alpha'
$$
with $B_\alpha'$ $\overset o\to\tau_i$-boxes.
Then
$$
\int_{B'} \phi^4_\tau \phi^4_{\tau_i} \sim \sum_\alpha\Big[\phi_{\tau_i}\Big|_{B_\alpha'} 
\Big]^4 \int_{B_\alpha'} \phi_\tau^4.\tag 3.1
$$

Note that $\overset o\to\tau_i$ is an $\big[\frac K\delta\times \frac K\delta\times \frac {K^2}{\delta^2}\big]$-box in direction $\xi_i$-normal at $\tau_i$.
Let $\xi$ be any normal for $\tau$.
Thus $\measuredangle (\xi, \xi_i)<\delta$ and $K\overset o\to\tau$ is contained in $\big[2\frac K\delta\times 2\frac K\delta\times 2\frac K{\delta^2}]$-box
in direction $\xi_i$.
It follows that $\overset o \to\tau_i$ may be partitioned in $K\overset o\to\tau$-boxes $B$ and hence by (2.4)
$$
\nint_{B_\alpha'}\phi^4_\tau \leq \max_B \nint_B \phi^4_\tau< C(K_1).\tag 3.2
$$
Substituting (3.2) in (3.1) gives
$$
C(K_1)\sum_\alpha\int_{B_\alpha'} \phi^4_{\tau_i}= C(K_1) \int_{B'} \phi^4_{\tau_i}< C(K_1)^2 |B'|.\tag 3.3
$$
Note also that in (2.2) $\mathcal L$ consists of at most $K$ $\frac \delta K$-discs.
Iteration of (2.1)-(2.3), where we terminate the process for (2.1) and continue for (2.2), gives a representation
$$
\align
&\qquad\qquad |Tf|\leq \\
& \qquad\qquad R^\ve\max_{1>\delta>\frac 1{\sqrt R}} \max_{\mathcal E_\delta}
\Big[\sum_{\tau\in\mathcal E_\delta} \Big(\phi_\tau|Tf_{\tau_1}|^{1/3} |Tf_{\tau_2}|^{1/3}
|Tf_{\tau_3}|^{1/3}\Big)^2\Big]^{\frac 12} \tag 3.4 \\ %\qquad\qquad (3.4)\\
&\qquad\qquad +\max_{\mathcal E_{\frac 1{\sqrt R}}}\Big[\sum_{\tau\in \mathcal E} (\phi_\tau|Tf_\tau|)^2\Big]^{\frac 12}
\tag 3.5 %\qquad\qquad\qquad\qquad\qquad (3.5)
\endalign
$$
where

\noindent
{\rm (3.6)}  \ $\mathcal E_\delta$ consists of at most $\frac 1\delta$ disjoint $\delta$-caps $\tau$

\noindent
{\rm (3.7)}  \ $\tau_1, \tau_2, \tau_3 \subset\tau$ are $\frac 1{K\delta}$-size and non-collinear

\noindent
{\rm (3.8)}  \ $\nint_B\phi^4_\tau < C(K_1)^{\frac{\log\frac 1\delta}{\log K}} < R^{\frac{\log C(K_1)}{\log K}} \ll R^\ve$
if $B$ is  a $\overset o\to\tau$-box.

Fix dyadic $1>\delta>\frac 1{\sqrt R}$ and consider
$$
\max_{\Cal E_\delta} \Big[\sum_{\tau\in\Cal E_\delta} (\phi_\tau|Tf_{\tau_1}|^{1/3} \, |Tf_{\tau_2}|^{1/3} \, |Tf_{\tau_3}|^{1/3})^2
\Big]^{\frac 12}\tag 3.9
$$
with $\Cal E_\delta$ and $\tau_1, \tau_2, \tau_3$ as above.

In what follows, we will make several estimates on (3.9) considering various norms.

\bigskip
\noindent
{\bf 4.} We assume $|f|\leq 1$.
By rescaling, for $\tau_1, \tau_2, \tau_3\subset\tau$ as in (3.7),
$$
\int_{B_R} |Tf_{\tau_1}|.|Tf_{\tau_2}|.|Tf_{\tau_3}|\leq \delta^2 \int_{B_{\delta R}} |Tg_{U_1}| \ |Tg_{U_2}| \ |Tg_{U_3}|.\tag 4.1
$$
with $|g|<1$ and $U_1, U_2, U_3\subset B_1$ of size $\sim\frac 1K$ and not collinear.

\medskip

Hence, from [BCT]
$$
\int_{B_{\delta R}} |Tg_{U_1}| \ |Tg_{U_2}| \ |Tg_{U_3}|\ll R^\ve\tag 4.2
$$
and
$$
\int_{B_R} |Tf_{\tau_1}|. |Tf_{\tau_2}|.|Tf_{\tau_3}|\ll \delta^2 R^\ve.\tag 4.3
$$
By (3.6) and H\"older
$$
\aligned
&\Big[\sum_{\tau\in\mathcal E_\delta} \big(\phi_\tau|Tf_{\tau_1}|^{\frac 13}. |Tf_{\tau_2}|^{\frac 13}.|Tf_{\tau_3}|^{\frac 13}\big)^2\Big]^
{\frac 12}\leq\\
&|\mathcal E_\delta|^{\frac 16}\Big[\sum_\tau \phi^3_\tau |Tf_{\tau_1}|. |T f_{\tau_2}|.|Tf_{\tau_3}|\Big]^{\frac 13}\leq\\
\endaligned
$$
$$
\delta^{-\frac 16} \Big[\sum_\tau \phi^3_\tau|Tf_{\tau_1}| \ |Tf_{\tau_2}| \ |Tf_{\tau_3}|\Big]^{1/3}
\tag 4.4
$$
where in (4.4) $\tau$ ranges over a partition in $\delta$-discs (note that (4.4) does not depend on $\mathcal E_\delta$ anymore).

We obtain
$$
\Vert(3.9)\Vert_{L^3(B_R)}\leq \delta^{-\frac 16}\Big[\sum_\tau\int\phi^3_\tau |Tf_{\tau_1}| \ |Tf_{\tau_2}| \
|Tf_{\tau_3}|\Big]^{1/3}.\tag 4.5
$$
Consider a partition of $B_R$ in $\overset o\to\tau$-boxes $B$.
Since $|Tf_{\tau_i}|$ are $\approx$ constant on $\overset o\to\tau_i$-boxes, hence on each $B$,
$$
\aligned
\int\phi^3_\tau|Tf_{\tau_1}| \ |Tf_{\tau_2}| \ |Tf_{\tau_3}|&\approx\sum_B(|Tf_{\tau_1}| \ |Tf_{\tau_2}| \ |Tf_{\tau_3}|)\Big|_B\Big(
\int_B\phi_\tau^3\Big)\\
&\approx \sum_B\Big[\int_B |Tf_{\tau_1}|.|Tf_{\tau_2}|.|Tf_{\tau_3}|\Big] \nint_B\phi^3_\tau\\
&\overset{(3.8)}\to\ll R^\ve\int_{B_R} |Tf_{\tau_1}|.|Tf_{\tau_2}| .|Tf_{\tau_3}|\\
\endaligned
$$
$$
\overset {(4.3)}\to < R^\ve \delta^2.\tag 4.6
$$
Therefore
$$
\Vert(3.9)\Vert_{L^3(B_R)} \ll R^\ve \delta^{-\frac 16}\tag 4.7
$$
which is our first bound.

\bigskip

\noindent
{\bf 5.} Take $3\leq p\leq 4$.

By H\"older again
$$
\aligned
(5.1)  = &\max_{\mathcal E_\delta} \Big[\sum_{\tau \in\mathcal E_\delta} (\phi_\tau |Tf_{\tau_1}|^{1/3} |Tf_{\tau_2}|^{1/3}
|Tf_{\tau_3}|^{1/3})^2\Big]^{1/2}\leq\\
& \Big(\frac 1\delta\Big)^{\frac 12-\frac 1p} \Big[\sum_\tau \phi_\tau^p(|Tf_{\tau_1}|.|Tf_{\tau_2}|.|Tf_{\tau_3}|)^{\frac p3}\Big]^{\frac 1p}
\endaligned
$$
implying
$$
\Vert(5.1)\Vert_p\leq \Big(\frac 1\delta\Big)^{\frac 12-\frac 1p} \Big[\sum_\tau\int_{B_R}\phi^p_\tau(|Tf_{\tau_1}| . |Tf_{\tau_2}|.|Tf_{\tau_3}|)^{p/3}
\Big]^{\frac 1p}.\tag 5.2
$$
As in (4.6)
$$
\aligned
&\int_{B_R}\phi^p_\tau(|Tf_{\tau_1}| .|Tf_{\tau_2}|.|Tf_{\tau_3}|)^{p/3}\leq\\
&\Big[\max_{B\, \overset o\to \tau-\text{box}} \ \nint_B\phi_\tau^p\Big] \ \Big[\int_{B_R}(|Tf_{\tau_1}| .|Tf_{\tau_2}| .|Tf_{\tau_3}|)^{p/3}\Big]\leq\\
&R^\ve\Big[\int_{B_R}|Tf_{\tau_1}|.|Tf_{\tau_2}|.|Tf_{\tau_3}|\Big]. \delta^{6(\frac p3-1)}\\
\endaligned
$$
$$
< R^\ve \delta^{2p-4}\tag 5.3
$$
by (3.8), (4.3) and since $\Vert Tf_{\tau_1}\Vert_\infty <\delta^2$.

Substituting (5.3) in (5.2) gives
$$
R^\ve \Big(\frac 1\delta\Big)^{\frac 12-\frac 1p}\Big(\frac 1\delta\Big)^{\frac 2p} \delta^{2-\frac 4p}=R^\ve \delta^{\frac 32 -\frac 5p}.\tag 5.4
$$
Hence
$$
\Vert(3.9)\Vert_{L^p(B_R)}\ll R^\ve \text { for } p\geq \frac{10}3 =p_0.\tag 5.5
$$
Returning to (5.1), let $0<\lambda<1$ be a parameter and denote
$$
g_{\tau}= |Tf_{\tau_1}|^{1/3}.|Tf_{\tau_2}|^{1/3}.|Tf_{\tau_3}|^{1/3} \text { and } \ g_{\tau, \lambda}= g_\tau 1_{[g_\tau \sim \lambda\delta^2]}  
\tag 5.6
$$
Then by (4.3)
$$
\int_{B_R}[g_{\tau, \lambda}]^p< (\lambda\delta^2)^{p-3} \int_{B_R}(g_{\tau, \lambda})^3 \ll R^\ve\lambda^{p-3} \delta^{2p-4}\tag 5.7
$$
and
$$
\Big\{\int_{B_R} \max_{\mathcal E_\delta}\Big[\sum_{\tau\in\mathcal E_\delta} (\phi_\tau g_{\tau, \lambda})^2 \Big]^{p_0/2}\Big\}^{1/p_0}
\ll R^\ve\lambda^{1-\frac 3{p_0}}=R^\ve \lambda^{\frac 1{10}}.\tag 5.8
$$
Let $1\leq \mu<\infty$ be another parameter and decompose each $\phi_\tau$ as
$$
\phi_\tau =\sum_{\mu \text { dyadic}}\phi_{\tau, \mu} \text { where }
$$
$$
\aligned
\phi_{\tau, \mu}=\phi_\tau 1_{[\phi_\tau\sim\mu]}\\
\phi_{\tau, 1} =\phi_\tau 1_{[\phi_\tau\leq 1]}\endaligned
\tag 5.9
$$
If $B$ is a $\overset o\to \tau$-box, (3.8) implies for $\mu>1$
$$
\nint_B \phi^{^{p_0}}_{\tau, \mu} \leq \mu^{-4+ p_0} \nint_B \phi^4_\tau \ll R^\ve \mu^{-2/3}.\tag 5.10
$$
Hence, instead of (5.8), we obtain
$$
\Big\{\int_{B_R}\max_{\mathcal E_\delta} \Big[\sum_{\tau \in\mathcal E_\delta} (\phi_{\tau, \mu} g_{\tau, \lambda})^2\Big]^{p_0/2} \Big\}^{1/p_0} \ll
R^\ve \lambda^{\frac 1{10}} .
\mu^{-\frac 15}.\tag 5.11
$$

Next, we perform a different type of estimate. Clearly
$$
\max_{\mathcal E_{\delta}} \Big[\sum_{\tau\in\mathcal E_\delta} (\phi_{\tau, \mu} g_{\tau, \lambda})^2 \Big]^{\frac 12} \leq \mu \Big(\sum_\tau g^2_{\tau,
\lambda}\Big)^{\frac 12}\tag 5.12
$$
with $\tau$ ranging over a partition in $\delta$-caps.

We apply the usual procedure to bound (5.12) by a Kakeya maximal function.

Writing
$$
|Tf_{\tau_i}|\lesssim |Tf_{\tau_i}| * (\delta^4 1_{\overset o\to \tau})
$$
we have
$$
\align
g_\tau(x) &\lesssim \int \Big\{\prod^3_{i=1} [|Tf_{\tau_i}|*(\delta^4 1_{\overset o\to\tau})]^{\frac 13}\Big\}(z) \
(\delta^4 1_{\overset o\to\tau}) (x-z)dz\\
&= \int \omega(z) (\delta^41_{\overset o\to\tau})(x-z) dz\tag 5.13
\endalign
$$
and
$$
g^2_{\tau, \lambda}(x) \lesssim \delta^4\int(\omega^2 1_{[\omega\gtrsim \lambda\delta^2]})(z) 1_{\overset o\to\tau}(x-z) dz.\tag 5.14
$$
Further
$$
\align
\int_{B_R} \omega^2 \ 1_{[\omega\gtrsim \lambda\delta^2]}&\lesssim \frac 1{\lambda\delta^2} \int\omega^3\\
&\lesssim \frac 1{\lambda\delta^2}\int\Big\{ \int\Big[\prod^3_{i=1} |Tf_{\tau_i}|(x-\tau_i)\Big]^{\frac 13} dx\Big\}
\Big[\prod^3_{i=1} (\delta^4 1_{\overset o \to\tau})(\tau_i)\Big] dz_1 dz_2 dz_3\\
\endalign
$$
$$
\ll R^\ve \lambda^{-1}.\tag 5.15
$$

Hence, from (5.14), (5.15), we obtain a representation
$$
g^2_{\tau, \lambda} \ll R^\ve \delta^4 \lambda^{-1} \int 1_{\overset o\to \tau}(\cdot -y) \Bbb P_\tau (dy).\tag 5.16
$$

>From  (5.16) and convexity
$$
\align
\Vert(5.12)\Vert_{L^{p_0}(B_R)} &\ll R^\ve \lambda^{-\frac 12}\mu\delta^2\Big\Vert\Big[\sum_\tau \ 1_{\overset o\to \tau}(x-y_\tau)\Big]^{\frac
12}\Big\Vert_{L^{p_0}{(B_R)}}\\
&= R^\ve \lambda^{-\frac 12}\mu\delta^2 \Big[\int\Big[\sum_\tau \ 1_{\overset o\to \tau} (x-y_\tau)]^{5/3}dx\Big]^{3/10}
\tag 5.17
\endalign
$$
for some choice of $\{y_\tau\}$-points in $\Bbb R^3$.

At this point we can invoke the $L^{5/2}$-bound for the $\Bbb R^3$-Kakeya maximal function.
In its dual formulation, we have
$$
\Big\Vert\sum_{v\in \frak S} \ 1_{T_v} \Big\Vert_{L^{5/3}}\leq \Big(\frac 1\kappa\Big)^{\frac 15+}\tag 5.18
$$
where $T$ is a translate of a tube of width $\kappa$ and length 1 in direction \hfill\break
$v\in\frak S~\subset~S_2$, where $\frak S$ consists of $\kappa$-separated points.

Rescaling by a factor $\delta^2$ and applying (5.18) with $\kappa=\delta$, it follows
$$
\Big\Vert\sum_\tau\ 1_{\overset o\to\tau}(\cdot-y_\tau)\Big\Vert_{L^{5/3}}\ll \delta^{-\frac{19}5-}.\tag 5.19
$$
Hence
$$
(5.17) \ll R^\ve \lambda^{-\frac 12}\mu\delta^{\frac 1{10}}.\tag 5.20
$$
which is our final estimate.

Summarizing (4.7), (5.11), (5.20), we have
$$
\Vert\max_{\Cal E_\delta}\Big[\sum_{\tau\in\Cal E_\delta} (\phi_\tau g_\tau)^2 \Big]^{\frac 12} 
\Big\Vert_{L^3 (B_R)} \ll R^\ve \delta^{-\frac 16}\tag 5.21
$$
and
$$
\align
\Big\Vert\max_{\Cal E_\delta}\Big[\sum_{\tau\in\Cal E_\delta} (\phi_{\tau, \mu} g_{\tau, \lambda})^2\Big]^{\frac 12}\Big\Vert_{L^{10/3}(B_R)} 
&\ll R^\ve \min \big(\lambda^{\frac 1{10}}\mu^{-\frac 15}, \lambda^{-\frac 12}\mu \delta^{\frac 1{10}}\big)\\
&\ll R^\ve \delta^{\frac 1{60}}\tag 5.22
\endalign
$$

Let
$$
q=\frac{33}{10}.
$$
Interpolating between (5.12), (5.22), it follows that
$$
\Vert(3.4) \Vert_{L^q(B_R)} \ll R^\ve.\tag 5.23
$$

\noindent
{\bf 6.} Remains to bound $\Vert(3.5)\Vert_q$.

Estimate
$$
\Big\Vert\Big[\sum_{\tau\in\Cal E} (\phi_\tau|Tf_\tau|)^2\Big]^{\frac 12}\Big\Vert_{L^3(B_R)} \leq (\sqrt R)^{\frac 16} \Big\{\sum_\tau\int_{B_R} \phi^3_\tau |Tf_\tau|^3 \Big\}^{\frac 13}
\tag 6.1
$$
where in the second sum, $\tau$ ranges over a full position in $\frac 1{\sqrt R}$-caps.

Since $|Tf_\tau|\lesssim \frac 1R$, (3.8) implies that
$$
(6.1) \ll R^{\frac 1{12}+\ve}.\tag 6.2
$$

On the other hand, using the decomposition (5.9), we obtain the following estimates on
$$
\Big\Vert\max_{\mathcal E_{\frac 1{\sqrt R}}}\Big[ \sum_{\tau\in\mathcal E}(\phi_{\tau, \mu}|Tf_\tau|)^2\Big]^{\frac 12}\Big\Vert_{L^{p_0}_{B_R}}.\tag 6.3
$$
Using (5.10), we get 
$$
\align
(6.3) &\leq (\sqrt R)^{\frac 12-\frac 1{p_o}}\Big(\sum_\tau\Vert \phi_{\tau, \mu} |Tf_\tau| \Vert^{p_0}_{L^{p_0}_{B_R}}\Big)^{\frac 1{p_0}}\\
&\ll (\sqrt R)^{\frac 12-\frac 1{p_0}+\ve} \mu^{-\frac 15} R^{\frac 1{p_0}}(\sqrt R)^{\frac 4{p_0}-2}\ll R^\ve \mu^{-1/5}.
\tag 6.4
\endalign
$$
Using the bound $\phi_{\tau, \mu} \lesssim \mu$ and the inequality
$$
|Tf_\tau|^2 \lesssim \frac 1{R^2} \int|Tf_\tau|^2 (y)\ 1_{\overset o\to \tau}(x-y)dy\tag 6.5
$$
and
$$
\int_{B_R} |Tf_\tau|^2 \lesssim 1\tag 6.6
$$
for $\tau\subset S_2$ a $\frac 1{\sqrt R}$-cap, we obtain similarly to (5.17)
$$
\align
(6.3)&\leq \mu\Big\Vert\Big(\sum_\tau |Tf_\tau|^2\Big)^{1/2} \Big\Vert_{L^{p_0}_{B_R}}\\
&\ll \frac\mu R \Big[\Big\Vert \sum_\tau 1_{\overset o\to \tau} (\cdot -y_\tau)\Big\Vert_{L^{5/3}_{B_R}}\Big]^{\frac 12}\\
&\overset {(5.19)}\to \ll \frac \mu R R^{\frac{19}{20}} = \mu R^{-\frac 1{20}+\ve}.
\tag 6.7
\endalign
$$
Hence, from (6.4), (6.7)
$$
(6.3) \ll R^\ve \min(\mu^{-\frac 15}, \mu R^{-\frac 1{20}}) \ll R^{-\frac 1{120} +\ve}. \tag 6.8
$$

Interpolation between (6.1), (6.8) implies
$$
\Vert(3.5)\Vert_{L^q(B_R)} \ll R^\ve.\tag 6.9
$$
Hence, we proved

\proclaim
{\bf Theorem 2$'$}
$$
\Vert Tf\Vert_{L^q(B_R)}\ll R^\ve \text { for } q\geq \frac{33}{10}, |f|\leq 1\tag 6.10
$$
\endproclaim
\noindent
(implying Theorem 2). 

\bigskip

\noindent
{\bf 7.} One can check how the preceding argument improves if one had the optimal Kakeya maximal function bound at disposal, thus
$$
\Vert\Cal M_\delta\Vert_{3\to 3} \ll \Big(\frac 1\delta\Big)^\ve\tag 7.1
$$
Recall (5.11)
$$
\Big\{\int_{B_R}\max_{\mathcal E_{\delta}}\Big[\sum_{\tau\in\mathcal E_\delta}(\phi_{\tau,\mu}g_{\tau, \lambda})^2
\Big]^{5/3}\Big\}^{3/10}  \ll R^{\ve}\lambda^{\frac 1{10}} \mu^{-\frac 15}.\tag 7.2
$$
Next, apply (5.17) with $p_0=3$
$$
\align
&\Big\Vert \max_{\mathcal E_\delta} \Big[\sum_{\tau\in\mathcal E_\delta} (\phi_{\tau, \mu} g_{\tau, \lambda})^2\Big]^{1/2}\Big\Vert_{L^3(B_R)}\ll\\
&R^\ve\lambda^{-\frac 12} \mu\frac 1R \Big[\int\Big[\sum_\tau \ 1_{\overset o\to\tau} 
(x-y_\tau)\Big]^{3/2} dx\Big]^{1/3} \ll R^\ve \lambda^{-\frac 12} \mu.
\tag 7.3
\endalign
$$
\medskip
For the (3.5) contribution, recall (6.3), (6.4)
$$
\Big\Vert\max_{\mathcal E_{\frac 1{\sqrt R}}} \Big[\sum_{\tau\in\mathcal E}(\phi_{\tau, \mu} |Tf_\tau|)^2\Big]^{\frac 12}\Big\Vert_{L^{10/3}_{(B_R)}}
\ll R^\ve\mu^{-1/5}\tag 7.4
$$
and using (6.5), (6.6), (7.1)
$$
\Vert \cdots \Vert_{L^3_{(B_R)}}\ll R^\ve\mu.\tag 7.5
$$
Interpolation between (7.2), (7.3) and (7.4), (7.5) gives
$$
\Vert Tf\Vert_{L^{q_1} (B_R)}\ll R^\ve  \ \text { for } \ q\geq \frac{36}{11} = 3,27 .... \text { and } \ |f|\leq 1.\tag 7.6
$$
This leads to an improved Theorem 2 with $3\frac 3{10}$ replaced by $\frac {36} {11}$.

\bigskip

\noindent
{\bf \S5. The Variable Coefficient Case}

We consider H\"ormander type oscillatory integral operators of the form
$$
(T_\lambda f)(x) =\int e^{i\lambda\psi(x, y)} f(y) dy\tag 5.1
$$
with real analytic phase function $\psi$ of the form
$$
\psi(x, y) =x_1y_1+\cdots+ x_{d-1} y_{d-1}+ x_d\big(\langle Ay, y\rangle + O(|y|^3)\big)+ O(|x|^2|y|^2)\tag 5.2
$$
and $\langle Ay, y\rangle$ a non-degenerate quadratic form.

Here $x$ (resp. $y$) are restricted to a neighborhood of $0\in\Bbb R^d$ (resp. $0\in\Bbb R^{d-1})$.
In order to bring (5.1) in the format considered earlier, rescale $x\to\frac x\lambda$ to obtain a phase function
$$
\phi(x, y) =x_1y_1+\cdots+ x_{d-1}y_{d-1}+x_d \big(\langle Ay, y\rangle +O(|y|^3)\big)+ \lambda\phi_\nu\Big(\frac x\lambda, y\Big)
\tag 5.3
$$
and $\phi_\nu$ at least quadratic in both $x, y$. Thus (5.1) becomes
$$
(Tf)(x) =\int e^{i\phi(x, y)} f(y) dy\tag 5.4
$$
with $x$ restricted to $|x|<o(\lambda)$.
This formulation appears as a perturbation of the restriction problem and preceding analysis can be generalized to
this setting.

First recall the [BCT] result in the variable coefficient case (see [BCT], Theorem 6.2 which treats the $d$-linear case,
but generalizes to lower levels of multi-linearity as formulated in [BCT], (40) for $\phi$ linear in $x$).

Thus let $1<k\leq d$ and
$$
(T_if)(x) =\int_{U_i} e^{i\phi_i(x, y)} f(y) dy\quad (1\leq i\leq k)\tag 5.5
$$
with $\phi_i$ as in (5.3).
We assume the transversality condition
$$
|Z_1 (x, y^{(1)}) \wedge\cdots\wedge Z_k(x, y^{(k)})|>c \text { for all $x$ and $y^{(i)}\in U_i$}\tag 5.6
$$
where
$$
Z(x, y)=\partial_{y_1} (\nabla_x \phi) \wedge\cdots\wedge\partial_{y_{d-1}} (\nabla_x\phi).\tag 5.7
$$
Then
$$
\Big\Vert\Big(\prod^k_{i=1} |T_if_i|\Big)^{\frac 1k} \Big\Vert_q \ll \lambda^\ve \Big(\prod^k_1\Vert f_i\Vert_2\Big)^{\frac 1k}\tag 5.8
$$
with $q=\frac{2k}{k-1}$ and $x$ restricted $|x|<o (|\lambda|)$.

Note that in the restriction problem, $Z(x, y)=Z(y)$ and (5.6) amounts to transversality of the normal vectors at the
corresponding hypersurface $S$ which is the graph of $\frac {\partial\phi}{\partial x_d}$.

It turns out that the $\lambda^\ve$-factor may be removed in (5.8) at the cost of increasing $q$ to $q_1>\frac {2k}{k-1}$.
Thus, as  proven in Lemma A3 in the Appendix, under the assumptions (5.5)-(5.7), one has
$$
\Big\Vert \Big(\prod^k_{i=1} |T_i f_i|\Big)^{\frac 1k} \Big\Vert_{q_1}  \leq C_{q_1} \Big(\prod^k_1 \Vert f_i\Vert_2\Big)^{1/k} \text { for }
q_1 >\frac {2k}{k-1}.\tag {5.8$'$}
$$
Using (5.8$'$) instead of (5.8) in \S2, \S3 to bound global multilinear contributions, will eliminate the $R^\ve$-factors (cf. \S3, (4.7) and
(4.9) for instance), without the need for an $\ve$-removal at the end (note that the $K^\ve$-factors  coming from a local application in \S3, (1.1) and
(2.2) are harmless).

\noindent
{\bf Remark.} We do not claim removal of the $\lambda^\ve$-factor in Theorem 6.2 from [BCT], but only in its present application to the
operators $T_i$ given by (5.5).

Returning to the analysis from \S2, \S3, also some adjustment is needed with respect to the parabolic rescaling argument that we discuss next.

Note that if we restrict $|y|<\frac 1K$ and rescale, letting $y=\frac{y'}K; x_1=Kx_1', \ldots, x_{d-1}=Kx_{d-1}'$
and $x_d=K^2 x_d'$, we obtain
$$
\int e^{i\phi'(x', y')} f\Big(\frac {y'}K\Big) dy' \text { where } \ |x_1'|, \ldots, |x_{d-1}'|<\frac \lambda K, |x_d'|<\frac \lambda{K^2}
\tag 5.9
$$
and
$$
\phi'(x', y')= x_1' y_1'+\cdots+ x_{d-1}' y_{d-1}' +x_d' \big(\langle A y', y'\rangle +\frac 1K O(|y'|^3)\big)
+\lambda\phi_\nu \Big(\frac{Kx_1'}{\lambda}, \ldots, \frac{Kx_{d-1}'}
{\lambda}, \frac {K^2x_d'}\lambda; \frac{y'}K\Big)
\tag 5.10
$$
with $x'$ subject to the restrictions (5.9).

Compared with (5.4), we see that one needs to consider the more general setting of operators
$$
(Tf)(x) =\int e^{i\phi(x, y)} f(y) dy \text { restricting } |x_1|, \ldots, |x_{d-1}|< R_1 \text { and } |x_d|< R\tag 5.11
$$
$(R\leq R_1)$, and
$$
\phi(x, y) =x_1 y_1 +\cdots + x_{d-1} y_{d-1}+x_d \big(\langle Ay, y\rangle+0(|y|^3\big)
 +R\phi_\nu \Big(\frac {x_1}{R_1}, \ldots, \frac {x_{d-1}}{R_1}, \frac {x_{d}}{R}; y\Big)
\tag 5.12
$$
(we use here that $\phi_\nu$ is at least quadratic in $y$).

It has to be shown that (5.8$'$) remains valid.
It turns out that the issue can be reduced to the $R=R_1$ case.
We give the details.
Let $q>\frac {2k}{k-1}$.

Partition the region
$$
Q=[|x_1|, \ldots, |x_{d-1}|< R_1] \times [|x_d|<R]=\bigcup_{s\leq \frac {R_1}{R}} Q_s
$$
in $R$-cubes and write
$$
\int_Q \Big(\prod^k_1 |T_i f_i|\Big)^{q/k} dx =\sum_s \int_{Q_s} \Big(\prod^k_1|T_i f_i|\Big)^{q/k} dx.\tag 5.13
$$
Partition the $y$-domain $\Omega\subset \Bbb R^{d-1}$ in cubes $\Omega_\alpha$ of size $\sim \frac 1R$ centered at $y_\alpha$ and write
$$
(T_if_i) (x) =\sum_\alpha^{(i)} e^{i\phi(x, y_\alpha)} \Big[\int_{\Omega_\alpha} f_i(y) e^{i[\phi(x, y)-\phi(x, y_\alpha)]} dy\Big].
$$
Restricting $x\in Q_s$, the factors [ \ ] are approximatively constant
$$
c_{i, \alpha}= \int_{\Omega_\alpha} f_i(y) e^{i[\varphi(\bar x, y)-\varphi(\bar x, y_\alpha)]} dy
$$
where $\bar x$ is the center of $Q_s$.
For $|z|<R$
$$
|T_if_i| (\bar x+z)\approx\Big|\sum_\alpha ^{(i)} e^{i\eta(z, y_\alpha)} e^{i\phi(\bar x, y_\alpha)} c_{i, \alpha}\Big|
$$
with $\eta(z, y_\alpha)= \phi(\bar x+z, y_\alpha)-\phi(\bar x, y_\alpha)$.
Hence, defining
$$
g_i(y) =c_{i, \alpha} e^{i\phi(\bar x, y_\alpha)} \text { for } y\in\Omega_\alpha
$$
we have 
$$
|T_if_i| (\bar x+z) \approx R^{d-1}\Big|\int e^{i\eta(z, y) } g_i(y) dy\Big|.
$$
>From (5.8$'$)
$$
\int_{B(0, R)} \Big[\prod^k_{1} |T_if_i| (\bar x+z)\Big]^{\frac qk}  \leq C R^{q(d-1)}\Big(\prod^k_{1} \Vert g_i\Vert_2\Big)^{\frac qk}
$$
$$
\aligned
& \leq CR^{\frac {q(d-1)}{2}}\Big[\prod^k_{1} \Big(\sum_\alpha^{(i)}|c_{i, \alpha}|^2\Big)^{\frac 12}\Big]^{\frac qk}
\\
&\leq C R^{\frac{q(d-1)}2} \Big\{\prod^k_{1} \Big[\sum_\alpha ^{(i)} \Big|\int_{\Omega_\alpha} f_i(y) 
e^{i\phi (\bar x, y)} dy\Big|^2\Big]^{\frac 12}\Big\}^{\frac qk}.
\endaligned
$$
Since $\bar x$ is any point in $Q_s$, we obtain
$$
R^{\frac{q(d-1)}2-d} \int_{Q_s} \Big\{\prod^k_{1}\Big[\sum_\alpha^{(i)} \Big|\int _{\Omega_\alpha} f_i(y) e^{i\phi(x,
y)}dy\Big|^2\Big]^{\frac 12}\Big\}^{\frac qk}.\tag 5.14
$$
Summing over $s$ gives
$$
\int_Q\Big[\prod^k_{1} |T_i f_i|\Big]^{\frac qk} < C R^{\frac{q(d-1)}2-d}\int_Q \prod^k_{1} \Big[\sum_\alpha^{(i)} \Big|\int_{\Omega_\alpha} f_i(y)
e^{i\phi(x, y)} dy\Big|^2\Big]^{\frac q{2k}}.\tag 5.15
$$

Note that
$$
\align
&\int_Q\Big|\int_{\Omega_\alpha} f_i(y) e^{i\phi(x, y)} dy\Big|^2\leq\\
&R.\max_{|x_d|<R} \int\Big| \int_{\Omega_\alpha} f_i(y) e^{i[x_1y_1+ \cdots + x_{d-1} y_{d-1}+R\phi_\nu
(\frac{x_1}{R_1},\cdots,  \frac {x_{d-1}}{R_1}, \frac{x_d} {R}; y)]} dy\Big|^2 dx_1 \cdots dx_{d-1}\\
&\lesssim R\int_{\Omega_\alpha}|f_i|^2\tag 5.16
\endalign
$$
using standard orthogonality considerations.

Also there is the trivial bound
$$
\align
\Big|\int_{\Omega_\alpha} f_i(y) e^{i\phi(x, y)}dy\Big|&\leq |\Omega_\alpha|^{\frac 12}(\int_{\Omega_\alpha}|f_i|^2\Big)^{\frac 12}\\
&\leq R^{-\frac {d-1} 2}\Big(\int_{\Omega_\alpha} |f_i|^2\Big)^{\frac 12}
\endalign
$$
implying
$$
\sum_\alpha\Big|\int_{\Omega_\alpha} f_i(y) e^{i\phi(x, y)}dy\Big|^2\lesssim R^{-(d-1)}\Vert f_i\Vert^2_2.\tag 5.17
$$

>From (5.15), (5.16), (5.17) and H\"older's inequality, it follows that 
$$
\align
(5.15) &\leq C R^{q \frac{d-1}2 - d} \Big\{ \prod^k_1(R^{1- \frac {d-1}2(q-2)} \Vert f_i\Vert_2^q)\Big\}^{1/k}
\\
&\leq C\Big(\prod^k_1\Vert f_i\Vert^q_2)^{1/k} \tag 5.18
\endalign
$$
as claimed.

We also observe that at suitable local scale, the phase function $\phi(x, y)$ given by (5.12) may be linearized in $x$,
reducing to the restriction setting.
Let $x=a+z\in B(a, \rho)$ and write
$$
\phi(x, y) =\phi(a, y) +\psi(z, y) +\Omega(z, y)\tag 5.19
$$
denoting
$$
\align
\psi (z, y)&= z_1y_1 +\cdots+ z_{d-1}  y_{d-1} + z_d \big(\langle Ay, y\rangle+ 0(|y|^3)\big)+\\
&\frac R{R_1} \Big\langle z', \nabla _{x'} \phi_\nu\Big(\frac {a'}{R_1}, \frac {a_d}{R}; y\Big)\Big\rangle
+ z_d \partial_{x_d} \phi_\nu \Big(\frac{a'}{R_1}, \frac {a_d}R; y\Big)\tag 5.20
\endalign
$$
with $x=(x', x_d)$ and where 
$$
 |\Omega (z, y)|= o(1) \text { provided $\rho =o(\sqrt R)$.}
\tag 5.21
$$
Since $\Omega$ does not oscillate on $B(a, \rho)$, it may be ignored in the phase function.

A suitable coordinate change in $y$ brings $\psi$ in the form
$$
\psi(z, y) =z_1 y_1+\cdots+ z_{d-1} y_{d-1} + z_d\big( \langle A'y, y\rangle + O(|y|^3)\big)
\tag 5.22
$$
with $A'$ a perturbation of $A$, hence $A'$ non-degenerate (and positive definite if $A$ is positive definite).

Using previous considerations, it is essentially straightforward to carry out the analysis from \S2, \S3 in the setting (5.11), (5.12),
assuming again that $A$ is positive definite and using (5.8$'$) to bound the global multilinear contributions.

Hence we obtain

\proclaim
{Theorem 3}
Consider the operator (5.1) with $\psi$ as in (5.2) and $A$ {positive definite}. Then
$$
\Vert T_\lambda f\Vert_{L^p_{\text {loc}}}\leq C_p \lambda^{-\frac dp}\Vert f\Vert_\infty\tag 5.23
$$
provided
$$
\aligned
p > 2\frac {4d+3}{4d-3} \quad & \text { if } d\equiv 0 (\mod 3)\\
p> \frac {2d+1}{d-1} &\text { if } d\equiv 1 (\mod 3)\\
p > \frac {4(d+1)}{2d-1} &\text { if } d\equiv2 (\mod  3).
\endaligned
$$
\endproclaim

In particular, for $d=3$, we obtain the condition $p > \frac{10}3$.
Interestingly, it turns out that this is the optimal exponent (as we will explain in the next section).

Without assuming $A$ positive definite, it is well-known that the condition
$$
p\geq \frac{2(d+1)}{d-1}\tag 5.24
$$
may be optimal range of validity for the inequality (5.19), when $d$ is {\it odd } (cf. [B]).

It was shown also in [B] that for $d$ even, there is some $p(d)<\frac{2(d+1)}{d-1}$ such that
$$
\Vert T_\lambda f\Vert_{L^p_{\text {loc}}} \lesssim \lambda^{-\frac dp} \Vert F\Vert_\infty. \tag 5.25
$$

The following statement makes this more precise

\proclaim
{Theorem 4}
Consider the operator (5.1) with $\psi$ as in (5.2) and $A$ non-degenerate.
For $d$ even, one has the inequality
$$
\Vert T_\lambda f\Vert_{L^p_{\text loc}} \leq C_p \lambda^{-\frac dp}\Vert f\Vert_\infty
\text { for } p > \frac {2(d+2)}d.\tag 5.26
$$
\endproclaim

(the exponent $\frac {2(d+2)}d$ was already known to be optimal).

\noindent
{\it Proof.} (sketch)

We consider the setting (5.11), (5.12).
Define the integer
$$
k=\frac d2+1.
$$
Thus the condition on the exponent $q$ in (5.8$'$) becomes  $q> \frac {2(d+2)}d$.

Following the procedure from \S2, \S3, we fix a large parameter $K$ and restrict $x$ to a $K$-ball $B_K=B(a, K)$.
Subdividing the $y$-domain $\Omega$ in balls $\Omega_\alpha$ of size $\frac 1K$ and considering the operators
$$
(T_\alpha f)(x) =\int_{\Omega_\alpha} e^{i\phi (x, y)} f(y)dy
$$
we consider the following two alternatives.
\bigskip

\noindent
{\bf Case 1.} On $B_K$, we may estimate
$$
|Tf| < C(K) \ |T_{\alpha_i} f|\tag 5.27
$$
for some $\alpha_1, \ldots, \alpha_k$ such that (5.6) holds for $y^{(1)} \in \Omega_{\alpha_1}, \ldots, y^{(k)}\in
\Omega_{\alpha_k}$ (with constant $c\sim \frac 1K$).
\bigskip

\noindent
{\bf Case 2.} Failure of Case 1.  This implies that on $B_K$
$$
|Tf|\lesssim \Big|\sum_{\alpha\in A} T_\alpha f\Big|+\max_\alpha |T_\alpha f|\tag 5.28
$$
where $\bigcup_{\alpha\in A} \Omega_\alpha$ is contained in an $\sim \frac 1K$-neighborhood of the $(k-2)$-manifold,
obtained by requiring $Z(a, y)$ given by (5.7) to belong to some $(k-1)$-dim linear space.

In particular,
$$
\#A\lesssim K^{k-2}.\tag 5.29
$$
In Case 1, write on $B_K$
$$
|Tf|\leq C(K) \sum_{\Sb \alpha_1, \ldots, \alpha_k\\ \text{(5.6) holds}\endSb} \Big(\prod^k_1|T_{\alpha_i} f|\Big)^{\frac
1k}.\tag 5.30
$$
The collected contribution may then be estimated using the $k$-linear bound and gives the estimate
$$
\ll C(K).\tag 5.31
$$
In Case 2, we proceed more crudely than in \S3 (note that lower dimensional restriction of the $y$-variable may
lead to degenerate phase functions if the quadratic form $\langle Ay, y\rangle$ is not assumed definite.)

>From (5.28)
$$
\align
\Big(\nint_{B_K} |Tf|^q\Big)^{\frac 1q} &\leq \Big( \nint_{B_K} \Big|
\sum_{\alpha\in A} T_\alpha f\Big|^q\Big)^{\frac 1q}
+ \Big(\sum_\alpha |T_\alpha f|^q\Big)^{\frac 1q}\\
&= (5.32)+(5.33)
\endalign
$$
Estimate
$$
\align
(5.32)^q &\leq \Big[\nint_{B_K} \Big|\sum_{\alpha\in A} T_\alpha f\Big|^2\Big] \ \Big[\sum_{\alpha\in A}
|T_\alpha f|\Big]^{q-2}\\
&\sim \Big[ \sum_{\alpha\in A} |T_\alpha f|^2 \Big] \ \Big[ \sum_{\alpha\in A} |T_\alpha f| \Big]^{q-2}
 \  \text { (using simple orthogonality)}\\
&< |A|^{1-\frac 2q+(q-2)(1-\frac 1q)} \sum_\alpha |T_\alpha f|^q.
\endalign
$$

Recalling (5.29)
$$
(5.32)\leq K^{(k-2) (1-\frac 2q)} \Big(\sum_\alpha \nint _{B_K} |T_\alpha f|^q\Big)^{1/q}\tag 5.34
$$
\big(that also captures (5.33)\big).

Thus the collected contribution over the $B_K$ is bounded by
$$
\align
&K^{(k-2)(1-\frac 2q)} \Big(\sum_\alpha \Vert T_\alpha f\Vert^q_q\Big)^{\frac 1q}\\
&\leq K^{(k-2)(1-\frac 2q)+\frac{d-1}q}\max_\alpha \Vert T_\alpha f\Vert_q.\tag 5.35
\endalign
$$
Rescaling gives the estimate
$$
<K^{(k-2)(1-\frac 2q)+\frac {d-1} q -(d-1)+\frac{d+1} q} Q^{(q)}_{\frac {R_1}K, \frac R{K^2}}
=K^{\frac {d+2}q -\frac d2} Q.\tag 5.36
$$
(denoting $Q^{(p)}_{R_1, R} $ a bound on $T:L^\infty\rightarrow L^p_{|x'|< R_1, |x_d|<R }$ given by (5.11)\big).

Since $q> \frac {2(d+2)}d$, this concludes the argument.

\bigskip
\noindent
{\bf \S6.  Some Examples}

We present in  this section an example for $n=3$ that will illustrate the optimality of the exponent $\frac {10}3$ in Theorem 3.
It will also explain the differences between the elliptic and hyperbolic cases.

Consider the following phase function
$$
\phi(x, y)= -x_1y_1-x_2y_2+\frac 12 x_3 y_1^2 +x_3^2 y_1y_2 +\frac 12(x_3+x_3^3) y_2^2.\tag 6.1
$$

First analyze the [BCT] transversality condition. Thus
$$
\nabla_x\phi =\Big(-y_1, -y_2, \frac 12(y_1^2+ y_2^2)+ 2x_3 y_1y_2+\frac 32 x^2_3 y_2^2\Big)
$$
$$
\left\{
\aligned &\partial _{y_1} \nabla_x\phi= (-1, 0, y_1+ 2x_3 y_2)\\
&\partial_{y_2} \nabla_x\phi = (0, -1, y_2+2x_3y_1+ 3x^2_3 y_2)
\endaligned
\right.
$$
$$
\aligned
Z(\Phi) (y, x)= \partial_{y_1} \nabla_x\phi\wedge\partial_{y_2}\nabla_x\phi&= (y_1+2x_3y_2, y_2+2x_3 y_1+ 3x_3^2 y_2, 1)\\
&=\Big(A\pmatrix y_1\\ y_2\endpmatrix, 1\Big)
\endaligned
$$
where $A=A_x=\pmatrix 1& 2x_3\\ 2x_3& 1+3x_3^2\endpmatrix$ is a perturbation of identity.

Concerning condition (40) in [BCT], if one fixes $x$ and restrict $y=(y_1, y_2)$ to non-collinear discs
$ U_1, U_2, U_3 \subset\mathbb R^2$, clearly
$$
\det\big(Z(\phi)(y^{(1)}, x), Z(\phi) (y^{(2)}, x), Z(\phi)(y^{(3)}, x)\big)\not= 0
$$
for $y^{(i)}\subset V_i$.
\bigskip

Next, consider the Kakeya type sets associated with (6.1).

$$
\left\{
\aligned \partial_{y_1} \phi&= -x_1+x_3 y_1+x_3^2 y_2\\
\partial _{y_2}\phi & = -x_2+ x_3^2 y_1+(x_3+x_3^3) y_2
\endaligned\right.
\tag 6.2
$$
and
$$
\Gamma_y \text { is parametrized by } \left\{\aligned x_1&= y_1x_3+y_2x_3^2\\
x_2& = y_1x_3^2+y_2(x_3+x_3^3).\endaligned\right.
\tag 6.3
$$
If we shift $\Gamma_y$ by $(y_2, 0, 0)$, the tubes
$$
\left\{\aligned x_1&= y_1x_3+ y_2x_3^2+y_2\\
x_2&= y_1 x_3^2+ y_2(x_3+x_3^3)\endaligned
\right.
\tag 6.4
$$
are contained in the surface
$$
S: x_1x_3=x_2.
$$
Thus one gets again 2D-compression, similar to the hyperbolic example
$$
\psi(x, y) =-x_1y_1-x_2y_2+2x_3 y_1y_2+x_3^2 y_2^2.\tag 6.5
$$
See also [Wi].
\bigskip

We try to exploit this compression as well as possible to make the oscillatory integral
$$
\int e^{i\lambda\phi(x, y)} \ f(y) dy\tag 6.6
$$
(with an appropriate $f$) large.

At this stage, there seems to be quite a difference between $(6.1)$ and $(6.5)$.
For $(6.5)$, just take
$$
f(y) =e^{iy_1^2}.\tag 6.7
$$
Then
$$
\int e^{i\lambda\psi(x, y)} f(y) dy= \int_{\text{loc}} e^{i\lambda[(y_1+x_3y_2)^2 -(x_1y_1+ x_2y_2)]}dy\tag 6.8
$$
and restricting $x$ to a $\frac 1\lambda$-neighborhood of $S$
$$
(6.8) \approx \int_{\text {loc}} e^{i\lambda[(y_1+x_3y_2)^2 -x_1(y_1+x_3y_2)]}dy.
$$
Setting $u=y_1+x_3y_2$, stationary phase implies
$$
|(6.8)|\sim\frac 1{\sqrt\lambda}.
$$
and hence
$$
\Vert(6.8)\Vert_{L_x^q} \sim \frac 1{\sqrt\lambda} \Big(\frac 1\lambda\Big)^{\frac 1q}\lesssim \Big(\frac 1\lambda\Big)^{\frac 3q}
 \text { for } q\geq 4.
$$

In the elliptic case, this type of construction seems impossible.

But one can make the following one, which will explain where the
condition $q\geq \frac{10}3$ comes from.

Instead of (6.7), take in (6.6)
$$
f(y) =\sum_{s < \sqrt\lambda} \sigma _s 1_{[\frac s{\sqrt\lambda}, \frac {s+c}{\sqrt\lambda}]} 
(y_2) e^{i\lambda \frac s{\sqrt\lambda} y_1}.\tag 6.9
$$
where $\sigma_s= \pm 1$ and $c>0$ is a small constant.

Hence
$$
(6.6) =\sum_{s < \sqrt\lambda} \sigma_s
\Big\{\int_{\frac s{\sqrt\lambda} < y_2 < \frac{s+c}{\sqrt\lambda}}
e^{i\lambda[\phi(x, y)+\frac s{\sqrt\lambda} y_1]} dy \Big\}.\tag 6.10
$$

Denoting $R$ the region
$$
R =\Big[x_3\sim 1\text { and } |x_2-x_1x_3|= o\Big(\frac 1{\sqrt\lambda}\Big)\Big]
\tag 6.11
$$
write
$$
\int_{\text{loc}} |(6.6)|^qdx\geq \int_R|(6.10)|^q dx.\tag 6.12
$$
Averaging the right side of (6.12) over signs $\sigma_s =\pm 1$, we obtain clearly
$$
\int_R\Big\{\sum_{s < \sqrt\lambda}\Big|\int_{\frac s{\sqrt\lambda} < y_2<\frac {s+c}{\sqrt\lambda}}
e^{i\lambda[\phi(x, y) +\frac s{\sqrt \lambda} y_1]}dy\Big|^2\Big\}^{\frac q2} dx.\tag 6.13
$$
Since
$$
\phi(x, y)=\frac 12 x_3\Big[\Big(y_1+x_3y_2 -\frac{x_1}{x_3}\Big)^2+\Big(y_2+\frac{x_1x_3-x_2}{x_3}\Big)^2\Big]-\frac 12
\Big[\frac{x_1^2}{x_3}+\frac{(x_1x_3-x_2)^2}{x_3}\Big]
$$
we have
$$
\phi(x, y) +\frac s{\sqrt\lambda} y_1=\frac 12 x_3 \Big[\Big(y_1+x_3y_2-\frac{x_1}{x_3}+\frac s{\sqrt\lambda}\frac 1{x_3}\Big)^2 +
\Big(y_2-\frac s{\sqrt\lambda}+\frac {x_1x_3-x_2}{x_3}\Big)^2\Big] + \eta(x, s)
\tag 6.14
$$
Therefore, from definition of $R$
$$
(6.13) \sim \int_R \Big\{\sum _{s < \sqrt \lambda}
\Big| \int_{\frac s{\sqrt\lambda}< y_2 <\frac{s+c}{\sqrt\lambda}}
e^{i\frac \lambda 2 x_3 (y_1 +x_3y_2 -\frac{x_1}{x_3} +\frac s{\sqrt\lambda} \frac 1{x_3})^2}
dy\Big|^2 \big\}^{\frac q2} dx.
\tag 6.15
$$

Stationary phase shows that for $|x_1| =o(x_3)$ and $s=o(\sqrt\lambda)$,
the inner integral in (6.15) is $O\Big(\frac 1\lambda\Big)$.

Hence
$$
(6.13) \sim \Big(\frac 1\lambda\Big)^{\frac {3q}4} |R| \sim \Big( \frac 1\lambda\Big)^{\frac {3q+2}4}
$$
by (6.11), and
$$
\Vert(6.6)\Vert_q \gtrsim \Big(\frac 1\lambda\Big)^{\frac 34 +\frac 1{2q}}.\tag 6.16
$$
Clearly (6.1) can only hold provided $q\geq \frac{10}3$.

\bigskip

\noindent
{\bf \S7. Curved Kakeya Estimates}

\noindent
{\bf 1.} Let's begin by describing curved Kakeya problems in $\mathbb{R}^n$.
We have a collection of tubes $T_i$.  Each tube $T_i$ is the $\delta$-neighborhood
of a curve $\Gamma_i$ in the unit ball in $\mathbb{R}^n$.  The goal of the curved
Kakeya problem is to assume some geometric information about the tubes $T_i$ and
use it to prove estimates for the $L^p$ norms of $\sum_i \chi_{T_i}$ and/or for
the volume of the union of tubes $\cup T_i$.  Either kind of estimate is a way of
measuring how much the tubes $T_i$ overlap.

Let $\delta>0$ be a small number.

We assume that each curve has $C^2$ norm $\lesssim 1$, and that each curve is an algebraic curve of
degree $ \lesssim 1$.  We assume that each curve is contained in the unit ball.  (I.e., $\Gamma_i$ is the restriction of an algebraic curve to
the unit ball.) (i.e. $\Gamma_i$ is the restriction of an algebraic curve to the unit ball.)

We define $T_i$ to be the $\delta$-neighborhood of $\Gamma_i$.  At each point $x \in T_i$, we can approximately define the tangent direction to
the tube $T_i$ at $x$.  Namely, pick any point $x' \in \Gamma_i \cap B(x, \delta)$ and define $v_i(x)$ to be the unit tangent vector to
$\Gamma_i$ at $x'$.  Since $\Gamma_i$ has $C^2$-norm $\lesssim 1$, choosing different points $x'$ in $B(x, \delta)$ will lead to an ambiguity of
size $\lesssim \delta$.  So the function $v_i(x)$ is well-defined up to $O(\delta)$ on the tube $T_i$.

\noindent
{\bf 2.} Assuming the $\Gamma_i$ algebraic, we prove the following slightly stronger version of the multilinear Kakeya estimate for curved tubes due to [BCT].
The  next statement deals with the 3-linear setting in $\Bbb R^4$  (for simplicity), but can be generalized to $k$-linear in $\Bbb R^n$.

\proclaim
{Theorem 6}

 Suppose $\Gamma_i$ are algebraic curves restricted to the unit 4-ball with degree $\lesssim 1$ and $C^2$ norm $\lesssim 1$.  
Let $T_i$ denote the $\delta$-neighborhood of $\Gamma_i$.  
Define approximate tangent vectors $v_i(x)$ for $x \in T_i$ as above. 
Suppose that the number
of tubes $T_i$ is $N$.  Then the following estimate holds:
$$ \int_{B^4} \Big[ \sum_{i=1}^N \chi_{T_i} \sum_{j=1}^N \chi_{T_j} \sum_{k=1}^N \chi_{T_k} \Big| v_i \wedge v_j \wedge v_k \Big|
\Big]^{1/2} \lesssim \delta^4 N^{3/2}.\tag 2.1 
$$
\endproclaim

Choosing the curves $\Gamma_i$ in the subspace $[e_1, e_2, e_3]$ implies immediately the same statement in $\Bbb R^3$ with 
$\delta^4$ replaced by $\delta^3$ in (2.1).

Since we may repeat tubes $T_i$, we obtain also the weighted version from Theorem~6.

The proof of the multilinear estimate follows the Dvir polynomial method, introduced
for problems over finite fields in [D].  The polynomial method was applied to multilinear
Kakeya problems in $\mathbb{R}^n$ in [G], and we will use results from there.

We will build an algebraic hypersurface $Z$ of controlled degree which is concentrated where
the tubes $T_i$ overlap heavily, and we will study the intersections between $Z$ and the curves
$\Gamma_i$.

Recall the definition of directed volume $V_{S} (v) := \int_S | v \cdot N | dvol_S$, where
$N$ denotes the normal vector to $S$.  We need a curved version of the cylinder estimate,
Lemma 2.1 in [G].

\proclaim
{Lemma 2.2}
If $Z$ is an algebraic surface in $\mathbb{R}^4$ of degree $D$, and if
$\Gamma_i$ is a curve of degree $d$, and if $Q_\alpha$ are disjoint cubes of side length
$\sim \delta$ which cover $T_i$, and if $x_\alpha$ is the center point of $Q_\alpha$,
then the following inequality holds:
$$ 
\sum_{\alpha} \delta^{-3} V_{Z \cap Q_\alpha} ( v_i(x_\alpha) ) \lesssim d D.\tag 2.3
 $$
\endproclaim

\noindent
{\it Proof.}  The idea of the proof is to interpret $\delta^{-3} V_{Z \cap Q_\alpha} (v_i (x_\alpha))$ in a
nice way: this quantity is roughly the average number of intersections of $Z \cap Q_\alpha$ with
a translation of $\Gamma_i$ by a random vector $v$ of length $\lesssim \delta$.  The total number
of intersections of $Z$ with (almost every) translate of $\Gamma_I$ is at most $d D$ by Bezout's theorem.

The errors caused by $v_i(x)$ varying by $\sim \delta$ as $x$ varies in
$Q_\alpha$ contribute about $\delta D$ per cube and so at most $D$ to the final answer.

In the paper [G], tubes had thickness 1.  
Our tubes have thickness $\delta$, so it's convenient to re-normalize
certain quantities.  If $Q \subset \mathbb{R}^4$ is a cube of side length $\delta$, then  
$$
V^{\text ren}_{Z \cap Q} (v) := \delta^{-3} V_{Z \cap Q} (v).\tag 2.4
$$
We recall the notion of `visibility' that plays a crucial role in [G].

The visibility of $Z \cap Q$ measures the directed volume of $Z \cap Q$ in various directions, and if
there is even one direction where $Z \cap Q$ has low directed volume, the visibility goes down a lot.  The
renormalized visibility has the following definition.
$$
Vis^{\text ren}[Z \cap Q] := Vol \Big( \{ v \text{ such that } |v| \le 1 \text{ and }
V^{\text ren}_{Z \cap Q}(v) \le 1 \} \Big) ^{-1}.\tag 2.5
$$

As in [G], one needs to introduce modified versions $\bar Vis$ and $\bar V$ of $Vis$ and $V$, obtained by a suitable
averaging over $Z$.
They have all good properties of the originals and moreover depend continuously on $Z$.
See [G] for details.

Next, we state a key result from [G] (see \S5, p14), in our renormalized setting.

\proclaim
{Lemma 2.6} Consider the standard $\delta$-lattice in $\Bbb R^4$.
Let $M$ be a function from the set of 4-cubes $Q$ in this lattice to $\Bbb Z_+\cup\{o\}$.
Then there is an algebraic hypersurface of degree $D$ such that
$$
\overline Vis^{\text {ren}}[Z\cap Q]\geq  M(Q)  \ \text { for all $Q$}\tag 2.7
$$
and
$$
D< C\Big[\sum_Q M(Q)\Big]^{1/4}.\tag 2.8
$$
\endproclaim

Let $Q_\alpha$ be a set of $\delta$-cubes that cover the unit 4-ball.  For each cube, define
$$
F(Q_\alpha) := \sum_{T_i, T_j, \text{ and } T_k \text{ intersect } Q_\alpha} |v_i \wedge v_j \wedge v_k|. $$

Here $v_i, v_j, v_k$ are evaluated at $x_\alpha$, the center of $Q_\alpha$.

\proclaim
{Lemma 2.9} The sum $\sum_\alpha \delta^4 F(Q_\alpha)^{1/2} \lesssim d^{3/2} \delta^4 N^{3/2}$.
\endproclaim

The sum on the left-hand side is very close to the integral over the 4-ball we want to estimate:
$$ 
\int_{B^4} \Big[ \sum_{i=1}^N \chi_{T_i} \sum_{j=1}^N \chi_{T_j} \sum_{k=1}^N \chi_{T_k} \left| v_i \wedge v_j \wedge v_k \right|
\Big]^{1/2} \sim \sum_\alpha  \delta^4 F(Q_\alpha)^{1/2}. 
\tag 2.10
$$

We compare our discrete sum and the integral below.  First we prove the lemma.

\noindent
{\it Proof.}
We construct a surface of degree $\lesssim D$ (for a large $D$) so that for all $\alpha$
$$
\bar Vis^{\text{ren}} [Z \cap Q_\alpha] \ge D^4 F(Q_\alpha)^{1/2} \Big[ \sum_\alpha F(Q_\alpha)^{1/2}\Big]^{-1}.\tag 2.11
 $$
(We can use any $D$, but we need $D$ big enough so that the RHS is at least 1 for all $\alpha$.)

The existence of $Z$ follows indeed from Lemma 2.6, taking for $M(Q_\alpha)$ the RHS of (2.11).

We show that
$$
D\Big[\sum_\alpha F(Q_\alpha)^{1/2}\Big]^{2/3} \lesssim d D N\tag 2.12
$$
which is equivalent with (2.9).
Write using (2.11).
$$ 
D \Big[ \sum F(Q_\alpha)^{1/2} \Big]^{2/3} \lesssim \sum F(Q_\alpha)^{1/3} \bar Vis^{ren} (Q_\alpha)^{1/3} D^{-1/3} \lesssim 
$$
$$ 
= \sum_\alpha \Big[ D^{-1} \bar Vis^{ren} (Q_\alpha) \sum_{T_i, T_j, T_k \text{ meet } Q_\alpha} |v_i \wedge v_j \wedge v_k (x_\alpha)|
\Big]^{1/3}.\tag 2.13 
$$

\noindent
{\bf Linear algebra lemma.}
 For any three vectors $v_i, v_j, v_k$, the following inequality holds
$$ 
\overline{Vis}^{ ren} [Z \cap Q_\alpha] |v_i \wedge v_j \wedge v_k| \lesssim D \bar V^{ren}_{Z \cap Q_\alpha}(v_i)
\bar V^{ren}_{Z \cap Q_\alpha}(v_j) \bar V^{ren}_{Z \cap Q_\alpha}(v_k)\tag 2.14
$$

{\it Proof.}
 We abbreviate $\bar V^{ren}_{Z \cap Q_\alpha}$ by $\bar V$ and $\overline{Vis}^{ren}$ by $\overline{Vis}$.

We use the following facts.  The function $\bar V$ maps $\mathbb{R}^4$ to $\mathbb{R}$.  It is non-negative.
It scales by the formula $\bar V ( \lambda v) = \lambda \bar V(v)$ for any $\lambda > 0$ and $v \in \mathbb{R}^4$.
It is convex.  And finally $|v| \le \bar V (v) \lesssim D |v|$ (where the lower bound is ensured by enlarging $Z$ with
$\sim \frac 1\delta $ hyperplanes.)

Now $\overline{Vis}$ is defined as $Vol \{ v \in B^4 | \bar V(v) \le 1 \}^{-1}$.  So we have to prove that
$$ 
Vol \{ v \in B^4 | \bar V(v) \le 1 \} \gtrsim |v_i \wedge v_j \wedge v_k | D^{-1} \bar V(v_i)^{-1} \bar V(v_j)^{-1} \bar V(v_k)^{-1}.
\tag 2.15
$$

Let $v_0$ be a unit vector perpendicular to the plane spanned by $v_i, v_j, v_k$.  Let $e_0 = v_0 / D$.  Then $\bar V (e_0) \le 1$.
Also, let $e_i := v_i / \bar V(v_i)$, so that $\bar V(e_i) = 1$.  Define $e_j, e_k$ similarly.  Since $\bar V(v) \ge |v|$,
it follows that $|e_i| \le 1$.  Since $\bar V$ is convex, $\bar V \le 1$ on the convex hull  of the eight points
$\pm e_0, \pm e_i, \pm e_j, \pm e_k$.  This convex hull lies in $B^4$.  Its volume is approximately
$|e_0 \wedge e_i \wedge e_j \wedge e_k|$.  Since $e_0$ is perpendicular to the other vectors, this wedge is
equal to  $|e_0| |e_i \wedge e_j \wedge e_k| = D^{-1} |v_i \wedge v_j \wedge v_k|  \bar V(v_i)^{-1} \bar V(v_j)^{-1} \bar V(v_k)^{-1}. $
proving (2.15).

\medskip

>From (2.14)
$$
(2.13)\lesssim \sum_\alpha \Big[ \sum_{T_i, T_j, T_k { meet } Q_\alpha}  \bar V^{ren}_{Z \cap Q_\alpha}(v_i)
\bar V^{ren}_{Z \cap Q_\alpha}(v_j) \bar V^{ren}_{Z \cap Q_\alpha}(v_k) \Big] ^{1/3} = 
$$
$$ 
= \sum_\alpha \sum_{T_i \text{ meets } Q_\alpha} \bar V^{ren}_{Z \cap Q_\alpha}(v_i) = \sum_{i=1}^N
\sum_{Q_\alpha \text{ meets } T_i}  \bar V^{ren}_{Z \cap Q_\alpha}(v_i).
$$

By the cylinder estimate, the last line is bounded $\lesssim N d D$ as required. 

This proves Lemma 2.9.

Finally, we return to the integral and show that the error in our discrete approximation is not too big:
$$ 
\int_{B^4} \Big[ \sum_{i=1}^N \chi_{T_i} \sum_{j=1}^N \chi_{T_j} \sum_{k=1}^N \chi_{T_k} \Big| v_i (x) \wedge v_j (x) \wedge v_k (x)
\Big| \Big]^{1/2} dx =$$

$$ = \sum_\alpha \int_{Q_\alpha} \Big[ \sum_{i=1}^N \chi_{T_i} \sum_{j=1}^N \chi_{T_j} \sum_{k=1}^N \chi_{T_k} \left| v_i \wedge v_j \wedge
v_k \right| \Big]^{1/2} dx $$

$$ \le \sum_\alpha \int_{Q_\alpha} \Big[ \sum_{T_i, T_j, T_k \text{ meet } Q_\alpha} |v_i (x) \wedge v_j(x) \wedge v_k(x)| \Big]^{1/2}
dx $$

$$ \le \sum_\alpha \int_{Q_\alpha} \Big[ \sum_{T_i, T_j, T_k \text{ meet } Q_\alpha} |v_i (x_\alpha) \wedge v_j(x_\alpha) \wedge
v_k(x_\alpha)| \Big]^{1/2} + \text { Error }\tag 2.16
 $$
where

$$ 
\text{Error }\lesssim \sum_\alpha \int_{Q_\alpha}  \Big[\sum_{i,j,k}  \chi_{\tilde T_i} \chi_{\tilde T_j} \chi_{\tilde T_k} |v_i \wedge v_j| \delta
\Big]^{1/2} \lesssim $$

$$ \lesssim \delta^{1/2} \Big( \int_{B^4} \sum_{i,j} \chi_{T_i} \chi_{T_j} |v_i \wedge v_j| dx \Big)^{1/2}
\Big( \int_{B^4} \sum \chi_{T_k} \Big)^{1/2} \sim $$

$$ N^{1/2} \delta^2 \Big( \int_{B^4} \sum_{i,j} \chi_{T_i} \chi_{T_j} |v_i \wedge v_j| dx \Big)^{1/2}.\tag 2.17
 $$

By Lemma 2.9, the first term in (2.16) is bounded by $C.\delta^4 d^{3/2}N^{3/2}$.

In (2.17) we encounter a 2-linear version of our original 3-linear integral.

This can be estimated by a much easier argument in the same spirit.

We show that
$$
\int_{B^4}\Cal X_{T_i} \Cal X_{T_j} |v_i\wedge v_j|< C\delta^4 d^2.\tag 2.18
$$
Hence $(2.17) < CN^{3/2}\delta ^4 d$ and this completes the proof of Theorem 8.

It remains to justify (2.18).
Thus

\proclaim
{Lemma 2.19}
Suppose that $\Gamma_1$ and $\Gamma_2$ are degree $d$ algebraic curves in $B^4$ and $C^2$
curves of norm $\lesssim 1$, and $T_i$ are $\delta$ tubes around $\Gamma_i$.

Then 
$$\int_{B^4} \chi_{T_1} \chi_{T_2} |v_1(x) \wedge v_2(x)| dx \lesssim d^2 \delta^4.
\tag 2.19
$$
\endproclaim

\noindent
{\it Proof.} (sketch) (This is an easier version of the 3-linear estimate (2.1)).

Cut the unit ball into $\delta$ cubes $Q_\alpha$.

Pick $D$ a large degree.  Choose $Z$ a degree $D$ hypersurface so that $\bar V^{ren}_{Z \cap Q_\alpha}(x) \ge |x|$
and
$$
\overline{Vis}^{ren} [ Z \cap Q_\alpha] \gtrsim D^4 |v_1 \wedge v_2 (x_\alpha) | \Big[ \sum_\alpha |v_1 \wedge
v_2 (x_\alpha) | \Big]^{-1}. 
\tag 2.20
$$

Now our integral is roughly
$$
 \delta^4 \sum_{Q_\alpha \subset T_1 \cap T_2} |v_1 \wedge v_2 (x_\alpha) |.\tag 2.21
$$

The error in this approximation is $\delta Vol (T_1 \cap T_2) \lesssim d \delta^4$ which is not larger
than the main term.

It suffices to prove 
$$
\sum_\alpha |v_1 \wedge v_2| \lesssim d^2.\tag 2.22
$$
Manipulating (2.20),  we see that

$$ \sum_\alpha |v_1 \wedge v_2| \lesssim D^{-4} [ \sum_\alpha \overline{Vis}^{ren}[Z \cap Q_\alpha]^{1/2}
|v_1 \wedge v_2|^{1/2} ]^2 \le $$

(by a linear algebra lemma like the one above)
$$
\lesssim D^{-2} [ \sum \bar V^{ren}(v_1)^{1/2} \bar V^{ren}(v_2)^{1/2} ] ^2 \le 
$$
$$
 \le D^{-2} \Big( \sum_\alpha \bar V^{ren}(v_1) \Big) \Big(\sum_\alpha \bar V^{ren}(v_2)\Big). 
$$

Now the first term in parentheses is essentially the average number of intersections between $Z$ and
$\Gamma_i$ after translating $\Gamma_i$ by a random vector of length $\lesssim \delta$, and so it
has size at most $d D$ by Bezout's  theorem.  (Compare the cylinder estimate above.)
The same applies to the second term.  So the whole
expression is bounded $\lesssim d^2$. 
\bigskip

\noindent
{\bf 3. Application to curved Kakeya sets}

Again we restrict ourselves to $n=4$ but the result generalize to even dimension $n$ (the exponent $\frac 32$ in Theorem 7 below
is then replaced by $1+\frac 2n$.)

Let the curves $\{\Gamma_i\}$ be as specified in the beginning of \S7.
We also make an `angle assumption' for pairs of curves, as follows.

The index set $\{ i \}$ is given a geometric structure.  For each curve $i$, we associate a point $y_i$
in $B^{n-1}(1)$.  We assume that the points $y_i$ are $\delta$-separated.  We make the following crucial geometric assumption.  If a point
$x$ lies in $T_i$ and in $T_j$, then the angle between $v_i(x)$ and $v_j(x)$ is $\gtrsim |y_i - y_j|$.  This assumption prevents too many
near-tangencies in the overlaps of the tubes.

\proclaim
{Theorem 7}
 Under the hypotheses above, for all $p > 3/2$,
$$ 
\Big\| \sum_i \chi_{T_i} \Big\|_p \lesssim \delta^{-3 + 4/p}.\tag 3.1
 $$
Hence, any curved Kakeya set in $\Bbb R^4$ (defined from algebraic curves of controlled degree and controlled $C^2$ norm) has Minkowski
dimension at least 3.\footnote"{$^{(*)}$}" { We will indicate later on in this section how to generalize this last claim to $C^\infty$-curves.}
\endproclaim

Examples (cf. [B2]) show that the statement in Theorem 7 is best possible.

The proof of Theorem 7 uses an inductive argument, where we assume that a good estimate holds
for a partial sums   $\sum_{{y_i}\in  \text{ small ball} } \chi_{T_i}$ and then we prove that a good estimate holds 
for a partial sum on $y_i$ in a larger ball.

\proclaim
{Theorem 7$'$}
 Let $T_i$ obey the hypotheses from Theorem 7 .  Suppose that $p > 3/2$.  Suppose that $\rho$ is a scale in the range $\delta \le
\rho \le 1$.  Let $B_\rho$ denote any ball of radius $\rho$ in $B^{3}(1)$.  Then the following estimate holds.
$$ 
\Big\| \sum_{y_i \in B_\rho} \chi_{T_i} \Big\|_p \lesssim \delta^{-3 + \frac{4}{p}} \rho^{3 - \frac{1}{p}}.\tag 3.2
 $$
\endproclaim

When $\rho = 1$, Theorem 7$'$ implies Theorem 7.  
When $\rho = \delta$, Theorem 7$'$ is trivial.  
We will prove Theorem 7$'$ by induction on $\rho$.
So we are allowed to assume that Theorem 7$'$ holds for all
$\bar\rho < \rho/2$.
In other words, we know
$$ 
\Big\| \sum_{y_i \in B_{\bar \rho}} \chi_{T_i} \Big\|_p \le \alpha \delta^{-3 + 4/p} \bar \rho^{3 - 1/p} .\tag 3.3
$$
In this equation $\alpha$ is a large constant that we will choose later.  Assuming (3.3), we will prove that the same estimate holds for balls of
radius $\rho$, {\it with the same constant $\alpha$}.  In other words, we will prove 
$$ 
\Big\| \sum_{y_i \in B_{\rho}} \chi_{T_i} \Big\|_p \le \alpha \delta^{-3 + 4/p} \rho^{3 - 1/p}. \tag 3.4 $$

Once we have proven (3.4), the inductive argument shows that Theorem 7$'$ holds for all $\rho$, and we are done.  
The idea of the proof is as follows.  We cover $B_\rho$ with smaller balls, and then write
$\sum_{y_i \in B_\rho}$ as a sum of contributions from the smaller balls.  
To bound the $L^p$ norm of this sum, we use a combination of two tools.  First,
(3.3) bounds the $L^p$ norms of the contributions from each smaller ball.  
By itself, this is not enough, but it shows that for (3.4) to fail, we
need to have points where many smaller balls are contributing.  The size of this effect is controlled by the multilinear estimate.

Let $K$ be a large constant to be determined later.  We cover $B_\rho$ by $K^3$ smaller balls, each of radius at most $10 \rho/K$.  We call each of these
smaller balls a ``clump".  Hence our set of tubes is divided into $\sim K^3$ clumps.

We divide $B^4$ into two regions, depending on how the tubes through $x$ are divided among the clumps.  We call a point $x \in B^4$ ``narrow'' if there exist $<
10^4 K$ clumps which contain half of the tubes through the point $x$.  We call $x$ ``broad" if it is not narrow.  Let $N \subset B^4$ be the set of narrow
points, and $N^c \subset B^4$ the set of broad points.

Our inductive hypothesis directly controls $\| \sum_{y_i \in B_\rho} \chi_{T_i} \|_{L^p(N)}$.

\proclaim
{Lemma 3.5}
Let $p > 3/2$.  Assuming (3.3), and assuming that $K = K(p)$ is sufficiently large, the following estimate holds:
$$
\int_{Narrow} \Big[ \sum_{y_i \in B_\rho} \chi_{T_i} \Big]^p dx \le (1/2) \alpha^p \delta^{4-3p} \rho^{3p-1}.\tag 3.6
$$

More explicitly, we say that $K$ is sufficiently large if $ [2 \cdot 10^7]^p K^{-2p + 3} < 1 /2$.  Notice that this condition
depends only on $p$.
\endproclaim

\noindent
{\it Proof.} Fix $x \in Narrow$.  We divide the sum $\sum_{y_i \in B_\rho} \chi_{T_i}(x)$ into clumps:
$$ 
\sum_{y_i \in B_\rho} \chi_{T_i}(x) \le \sum_{j=1}^{K^3} \Big[\sum_{y_i \in clump(j)} \chi_{T_i}(x) \Big].
\tag 3.7
$$

Now since $x$ is narrow, the sum on the right-hand side is controlled by the sum from only $10^4 K$ clumps.  In other words, we can pick a set $C(x)$ of at
most $10^4 K$ clumps so that
$$ 
(3.7) \le 2 \sum_{j \in C(x)} \big[  \sum_{y_i \in { clump }(j)} \chi_{T_i}(x) \Big].\tag 3.8
 $$

Now by Holder's inequality, this last sum is dominated by
$$ 
(3.8) \le 2  \Big[ \sum_{j \in C(x)} \Big( \sum_{y_i \in  {clump }(j)} \chi_{T_i}(x) \Big)^p \Big]^{1/p} [10^4 K]^\frac{p-1}{p} . $$

Putting together the string of inequalities we just proved, we see that for each $x \in Narrow$,
$$  
\Big[ \sum_{y_i \in B_\rho} \chi_{T_i}(x)\Big]^p \le 2^p  [10^4 K ]^{p-1} \sum_{j =1}^{K^3} \Big( \sum_{y_i \in {clump} (j)} \chi_{T_i}(x) \Big)^p. 
$$

Now integrating over the narrow set, we get
$$
\int_{Narrow}  \Big[ \sum_{y_i \in B_\rho} \chi_{T_i}(x)\Big]^p dx \le 2^p [10^4 K]^{p-1} \sum_{j=1}^{K^3}
\int_{B^4} \Big(\sum_{y_i \in {clump}(j)} \chi_{T_i}(x) \Big)^p dx. 
\tag 3.9
$$

But by induction (3.3), the integral involving each smaller clump in (3.9) is controlled
$$
\int_{B^4} \Big(\sum_{y_i \in clump (j)} \chi_{T_i}(x) \Big)^p dx \le \alpha^p \delta^{4 - 3p} (10 \rho / K )^{3p-1}. $$

Plugging this estimate into (3.8), we get
$$
\int_{Narrow} \Big [ \sum_{y_i \in B_\rho} \chi_{T_i}(x)\Big]^p dx \le 2^p [10^4 K]^{p-1} K^3
\alpha^p \delta^{4 - 3p} (10 \rho / K )^{3p-1}.\tag 3.10
$$

Grouping terms in the right-hand side, we get
$$ 
\le [2 \cdot 10^4 \cdot 10^3]^p K^{-2p+3} \alpha^p \delta^{4-3p} \rho^{3p-1}. $$

We choose $K = K(p)$ sufficiently large so that
$$ 
[2 \cdot 10^7]^p K^{-2p + 3} < 1 /2 .$$

Since $p > 3/2$, we can choose $K$ sufficiently large to make this inequality hold.
This proves Lemma 3.5.

At this point we fix $K = K(p)$.

Next we have to control the contribution from the broad points in $B^4$.  We do this using the multilinear estimate.

\proclaim
{Lemma 3.11}
Let $Broad \subset B^4$ denote the set of broad points.
$$ \int_{Broad} \Big| \sum_{y_i \in B_\rho} \chi_{T_i}\Big|^{\frac{3}{2}} \le C(K) \delta^{-1/2} \rho^{7/2}.\tag 3.12 $$
\endproclaim

\noindent
{\it Proof.} Let $x \in B^4$ be a broad point.  The broadness of $x$ leads to the following estimate:
$$ 
\Big| \sum_{y_i \in B_\rho} \chi_{T_i}(x)\Big|^3 
\le  \rho^{-2} C(K) \sum_{y_i \in B_\rho} \chi_{T_i}(x) \sum_{y_j \in B_\rho} \chi_{T_j}(x)  \sum_{y_k \in
B_\rho} \chi_{T_k}(x) | v_i(x) \wedge v_j(x) \wedge v_k(x) |  .\tag 3.13
$$

This holds because most triples of tubes through a broad point lie in clumps that fail to be coplanar, and so we have $|v_i (x) \wedge v_j(x) \wedge v_k(x)|
\ge \rho^{2} / C(K)$.

Taking the square root of (3.13) and integrating, we get
$$ 
\int_{Broad}  \Big| \sum_{y_i \in B_\rho} \chi_{T_i}\Big|^{\frac{3}{2}} \le 
$$
$$ 
\le C(K) \rho^{-1}  \int_{B^4} \Big[ \sum_{y_i \in B_\rho} \chi_{T_i}(x) \sum_{y_j \in B_\rho} \chi_{T_j}(x)  \sum_{y_k \in B_\rho} \chi_{T_k}(x) |
v_i(x) \wedge v_j(x) \wedge v_k(x) |  \Big]^{1/2} dx.
\tag 3.14
$$

But the right-hand side is controlled by the Multilinear Estimate.  The number of points $y_i \in B_\rho$ is $\le 100 [ \rho / \delta ]^3$.  
According to Theorem 6, the right-hand side is bounded above by
$$
(3.14) \lesssim_K \rho^{-1} \delta^4 [\rho / \delta]^{9/2}  = \delta^{-1/2} \rho^{7/2}
$$ 
proving Lemma 3.11.

The estimate in Lemma 3.11 controls the $L^{3/2}$ norm of $\sum \chi_{T_i}$ on the broad set.  There is
an obvious estimate for the $L^\infty$ norm, and by combining them we can estimate the $L^p$ norm
for our choice of $p > 3/2$.

We clearly have the $L^\infty$ bound
$$ 
\sup_x \Big|\sum_{y_i \in B_\rho} \chi_{T_i} \Big| \lesssim \rho^3 \delta^{-3} .\tag 3.15
$$

Since our $p > 3/2$, we see that
$$
\int_{Broad} \Big|\sum_{y_i \in B_\rho} \chi_{T_i}\Big|^p dx 
\lesssim [\rho^3 \delta^{-3}]^{p-3/2} \int_{Broad} \Big|\sum_{y_i \in B_\rho} \chi_{T_i}\Big|^{3/2} dx. $$

Applying Lemma 3.11 to bound the last integral, we see that
$$ 
\int_{Broad} \Big|\sum_{y_i \in B_\rho} \chi_{T_i}\Big|^p dx \le C(K) \rho^{3p-1} \delta^{4-3p}.\tag 3.16
 $$

Now we choose $\alpha$ large enough that $C(K) \le (1/2) \alpha^p$.  (So $\alpha$ depends on $K$ and $p$.)  Now we know that
$$ 
\int_{Broad} \Big|\sum_{y_i \in B_\rho} \chi_{T_i}\Big|^p dx \le (1/2) \alpha^p \rho^{3p-1} \delta^{4-3p}.\tag 3.17
$$
and (3.6), (3.17) imply (3.4).

This concludes the proof of Theorem 7$'$ and hence Theorem 7.

\bigskip
\noindent
{\bf 4. Estimates for $C^k$ curves}

We can prove estimates for $C^k$ curved Kakeya sets by approximating
the $C^k$ curves using polynomials.  This idea was suggested to us 
by Alex Nabutovsky.  He referred us to Jackson's theorem in approximation theory
and related results.

The results in this section look far from optimal, but we wanted to show
that something can be done for non-algebraic curves as well with these methods. 

\proclaim
{Jackson type theorem}
If $f: [0,1] \rightarrow \mathbb{R}$ has $C^k$ norm 1, then
we can approximate $f$ by a degree $d$ polynomial $P$ so that
$$ 
|f(x) - P(x) | \lesssim d^{-k} \text{   for all } x \in [0,1].\tag 4.1
$$
\endproclaim

In particular, we may approximate a $C^k$ curve $\Gamma_i$ by a degree $d$
algebraic curve with the same $\delta$-tube and with $d \lesssim \delta^{-1/k}$.

\noindent
{\bf Remark.} This algebraic curve will be just the graph of a degree $d$ polynomial.
There are many more algebraic curves and so one may hope for a better estimate,
but it would take some more sophisticated approximation theory.

Tracking the dependence on degree in Theorem 7, the following estimate is gotten.

\proclaim
{Theorem 7$''$}
Under the hypotheses in section 3, for all $p > 3/2$,
$$ 
\Big\Vert \sum_i \chi_{T_i} \Big\Vert _p \lesssim d^{\frac{3}{2p}} \delta^{-3 + 4/p}.
\tag 4.2 
$$
\endproclaim

Hence we get the following estimate for $C^k$ curves $\Gamma_i$ with
$k \ge 2$ obeying the angle condition:

\proclaim
{Theorem 8} Under the hypotheses above, for all $p > 3/2$,
$$
 \Big\| \sum_i \chi_{T_i} \Big\|_p \lesssim_k \delta^{-3 + \frac{4}{p} - \frac{3}{2 p k}}.
\tag 4.3
 $$
\endproclaim

In particular, for $C^\infty$ curves we have essentially the same estimate
that we had for algebraic curves.

An immediate consequence of Theorem 8 is the following result on the Minkowski dimension of curved Kakeya sets.

\proclaim
{Theorem 9}
Any curved Kakeya set in 4D associated to $C^\infty$-curves obeying the angle condition, has Minkowski dimension at least 3.
\endproclaim

The method described in \S7 can be generalized to higher dimension.
In particular, for $n$ even, smooth curved Kakeya sets in $\Bbb R^n$ have Minkowski dimension at least $\frac n2+1$.
This statement, which in some sense is the companion to Theorem 4, is the sharp version of a phenomenon first observed in [B2].
Note that for $n$ odd, (algebraic) curved Kakeya sets may have Minkowski dimension $\frac{n+1} 2$ (cf. [B2]).

\bigskip

\noindent
{\bf \S8. Further Comments}

It is not quite clear at this point what is the exact potential of the method introduced in this paper (when the optimal result is not attained) and we
have not tried to push the techniques to their limit.  In particular, further improvements in Theorem 2 are not out of question and one could also
explore if the more refined strategy used to obtain Theorem 2 in 3D has a higher dimensional counterpart (possibly improving upon Theorem 1).

Returning to inequality (5.8$'$) in \S5, we present next an alternative proof for $n=3$ of the following statement (which suffices for the application
to Theorem 3 when $n=3$).

\proclaim
{Proposition 8.1}
Under the transversality assumption (5.6), (5.7) from \S5 one has the 3-linear estimate in $3D$
$$
\Big\Vert\prod^3_{i=1} (T_\lambda^{(i)} f_i)\Big\Vert_{L^{q/3}} < \lambda^{-\frac 9q} \prod^3_{i=1} \Vert f_i\Vert_2\text { for }
q>\frac {10}3\tag 8.2
$$
where the operators $T_\lambda^{(i)}$ are given by (5.1), (5.2) with positive definite quadratic form and the phase functions are assumed algebraic of
bounded degree.
\endproclaim

Proposition 8.1 is weaker then (5.8$'$) in \S5, but may be obtained directly without the need for an $\ve$-removal lemma;
hence this argument may have some interest.

Returning to the argument in [BCT] (which is similar to the one in [B1]) there are basically two steps, that will be suitably modified.

\noindent
{\bf 1.} The first step in the approach involves the `intermediate scale' $|x|<\frac 1{\sqrt\lambda}$.
At this scale, as explained in (5.19)-(5.23) from \S5, the problem may be linearized in $x$.
This allows to derive a trilinear bound from the bilinear $2\times 2\to\frac q2$ estimate for $q>\frac {2(d+1)}d =\frac {10} 3$ due to [T1] in
the restriction theory rather than relying on a bootstrap.
We point out that the linear result from [T1] for the paraboloid and, more generally, smooth hypersurfaces with positive definite second fundamental
form, may fail without this last hypothesis (for instance for a hyperbolic paraboloid, cf. [V]), if no additional assumptions.

\noindent
{\bf 2.} 
At the second stage of the argument, the issue is the 3-linear Kakeya estimate (in the curved case), which is Proposition 6.8 in [BCT].
Here another factor $\lambda^\ve$ enters in their argument.
However, Theorem 6 of the paper may be used, since it immediately implies (by lowering the dimension from $\Bbb R^4$ to $\Bbb R^3$).

\proclaim
{Proposition 8.3}
Denoting $\{\tau_i \}$ $\delta$-neighborhoods of a family $\{\Gamma_i\}$ of smooth algebraic curves of degree $\lesssim 1$ in $B(0, 1)\subset\Bbb R^3$
and $v_i$ the tangent vector at a given point $p\in\Gamma_i$, one has
$$
\int\Big[\sum_{i, j, k} \lambda_i \mu_j \eta_k\Cal X_{\tau_i\wedge\tau_i \wedge\tau_k}|v_i\wedge v_j\wedge v_k|]^{1/2}< C\delta^3
\Big(\sum|\lambda_i|\Big)^{1/2} \Big(\sum |\mu_j|\Big)^{1/2} \Big(\sum|\eta_k|\Big)^{1/2}
$$
\endproclaim

\noindent
{\bf Proof of Proposition 8.1}

Rescaling $x\to \frac x\lambda$, we obtain the phase function
$$
\phi(x, y)= \lambda\phi\Big(\frac x\lambda, y \Big) \text { where } |x|=o(\lambda).
$$
Partition the $y$-domain $\Omega$ in boxes $\Omega_\alpha$ of size $\frac 1{\sqrt\lambda}$ centered at points $y_\alpha$.
Write for $y\in \Omega_\alpha$
$$
\phi(x, y)=\phi(x, y_\alpha)+\langle\nabla _y\phi (x, y_\alpha)\rangle _O(\lambda|y-y_\alpha|^2)
$$
where the last term may be dropped.

$$
T_\alpha f(x) =\int_{\Omega_\alpha} e^{i\langle \nabla_y\phi(x, y_\alpha), y-y_\alpha\rangle} f(y) dy
\tag 8.4
$$
and write
$$
Tf(x) =\sum_\alpha  e^{i\phi(x, y_\alpha)} (T_\alpha f)(x).\tag 8.5
$$
Next, introduce a variable $z\in B(0, \sqrt\lambda)$, writing
$$
Tf(x+z)\sim \sum_\alpha e^{i\phi(x+z, y_\alpha)} (T_\alpha f)(x).\tag 8.6
$$
Returning to (8.2), write
$$
\int_{B(0, \lambda)} \Big[\prod^3_{i=1} |T^{(i)} f_i|\Big]^{\frac q3} \sim \lambda^{-\frac 32}\int_{B(0, \lambda)} \Big\Vert\prod^3_{i=1}
(T^{(i)} f_i)(x+z)
\Big\Vert^{q/3}_{L^{q/3}(|z|<\sqrt \lambda)}dx\tag 8.7
$$
with $T^{(i)} f_i(x+z)$ replaced by (8.6).

Estimate $\big\Vert\prod^3_{i=1}\big\Vert_{\frac q3} \leq\big\Vert\prod_{i=1, 2}\big\Vert^{\frac 12}_{\frac q2} \big\Vert\prod_{i=2, 3}\big\Vert
^{\frac 12}_{\frac q2} \big\Vert\prod_{i=3, 1} \big\Vert^{\frac 12}_{\frac q2}$.

Denoting
$$
\eta(z, y) =\phi(x+z, y)-\phi(x, y) \qquad (x \text { fixed})
$$
we bound
$$
\int_{B(0, \sqrt \lambda)} \Big|\sum_\alpha e^{i\eta (z, y_\alpha)} (T^{(1)}_\alpha f_1)(x)\Big|^{\frac q2} \ \Big|\sum_\beta
e^{i\eta(z, y_\beta)} (T_\beta^{(2)} f_2)(x)\Big|^{\frac q2} dz.\tag 8.8
$$
Define functions $g_1, g_2$ by
$$
g_1(y) = e^{i\eta(z, y_\alpha)}(T_\alpha^{(1)} f_1)(x) \text { for } y\in\Omega_\alpha\tag 8.9
$$
and similarly for $g_2$.

Clearly
$$
(8.8) \sim\lambda^{q}\int_{B(0, \sqrt \lambda)} \Big|\int e^{i\eta (z, y)}g_1(y)dy\Big|^{\frac q2} \ \Big|\int e^{i\eta(z, y)}
g_2(y)dy\Big|^{\frac q2} dz.\tag 8.10
$$
Since, following (5.19)-(5.22) in \S5, $\eta$ has the form
$$
\eta (z, y) =z_1y_1+z_2y_2+z_3\big(\langle Ay,y \rangle +O(|y|^3)\big) +O(|z|\frac{|x|}\lambda |y|^2)+ O\Big(\frac{|z|^2}\lambda |y|^2\Big)\tag 8.11
$$
the last term in (8.11) may be dropped for $|z|<\sqrt\lambda$.
Hence $\eta(z, y)$ may be viewed as linear in $z$, of the form
$$
z_1y_1+ z_2y_2 +z_3\langle A'y, y\rangle +O(|z| \ |y|^3)\tag 8.12
$$
with $A'$ positive definite.

Applying the bilinear $2\times 2 \to\frac q2$ bound from [T1], it follows that
$$
\align
(8.10) &\lesssim \lambda^{q} \Vert g_1\Vert_2^{\frac q2} \Vert g_2\Vert_2^{\frac q2}\\
&\sim \lambda^{q/2} \Big[\sum_\alpha |(T_\alpha^{(1)} f_1)(x)|^2\Big]^{\frac q4}
\Big[\sum_\beta|(T_\beta^{(2)} f_2)(x)|^2\Big]^{\frac q4}.\tag 8.13
\endalign
$$

>From (8.13), the following bound on (8.7) is obtained
$$
\lambda^{\frac 12 (q-3)} \int_{B(0, \lambda)} \Big\{\prod^3_{i=1} \Big[\sum_\alpha|(T_\alpha^{(i)} f_i)(x)|^2\Big]^{\frac q6}\Big\}dx.
\tag 8.14
$$
The next step is to capture the factors in (8.14) by curved Kakeya maximal functions.
>From definition of $T_\alpha$
$$
|T_\alpha f|^2(x)=|\hat f_\alpha|^2 \big(\nabla_y \phi(x, y_\alpha)\big)\text { where }  f_\alpha = f|_{\Omega_\alpha}.\tag 8.15
$$
Let $b$ be a standard bumpfunction on $\mathbb R^{d-1}$.
Then $|\hat{f_\alpha}|^2$ may be recovered by an average of translates $b(\frac {\xi -\cdot}{\sqrt\lambda})$ with averaging weight $\lambda^{-1}
\Vert f_\alpha\Vert_2^2$.

Denoting
$$
c_\alpha^{(i)} =\Vert f_{i, \alpha}\Vert^2_2 \qquad (i=1, 2, 3)
$$
satisfying
$$
\sum_\alpha c_\alpha^{(i)} =\Vert f_i\Vert^2_2\tag 8.16
$$
we obtain therefore
$$
\sum_\alpha |(T_\alpha^{(i)} f_i)(x)|^2 \lesssim \lambda^{-1} \sum^{(i)}_{\alpha,\nu} b \big(\lambda^{-\frac 12}
(\nabla_y \phi(x, y_\alpha)-\xi_{\alpha, \nu})\big). c^{(i)}_{\alpha, \nu}
\tag 8.17
$$
where $c^{(i)}_{\alpha, \nu}> 0, \sum_\nu c^{(i)}_{\alpha, \nu}= c_\alpha^{(i)}$ and
$$
\sum_{\alpha, \nu} c^{(i)}_{\alpha, \nu} =\Vert f_i\Vert_2^2.\tag 8.18
$$
Substituting (8.17) in (8.14), one gets
$$
\lambda^{-\frac 32}\int_{B(0, \lambda)} \Big\{\prod^3_{i=1} \Big[\sum_{\alpha,\nu}^{(i)} b\big(\lambda^{-\frac 12}(\nabla_y\phi(x,
y_\alpha)-\xi_{\alpha, \nu})\big) c^{(i)}_{\alpha, \nu}\Big]^{\frac q6}\Big\}dx\\
$$
$$
\qquad =\lambda^{\frac 32}\int_{B(0, 1)}\Big\{\prod^3_{i=1} \Big[\sum_{\alpha, \nu}^{(i)} c^{(i)}_{\alpha, \nu} b\big(\lambda^{-\frac 12}
(\nabla_y\phi(\lambda x', y_\alpha) - \xi_{\alpha, \nu})\big)\Big]^{\frac q6}\Big\} dx'
\tag 8.19
$$
We may now apply Proposition 8.3.
In the present trilinear setting, $|v_i\wedge v_j\wedge v_k|>c$ and hence
$$
\Big\Vert\prod^3_{i=1} \Big[\sum_s \lambda_s^{(i)} \Cal X_{{\tau_s}^{(i)}}\Big]\Big\Vert _{L^{1/2}}\leq
c\delta^6\prod^3_{i=1} \Big[\sum_s|\lambda_s^{(i)}|\Big]\tag 8.20
$$
where $\delta=\frac 12$ and the tubes $\tau$ of the form
$$
\lambda^{-1}|\nabla_y \phi (\lambda x, y)-\xi|<\lambda^{-1/2}.\tag 8.21
$$
Interpolation of (8.20) with the obvious $L^\infty$-bound gives, for $r\geq \frac 12$
$$
\Big\Vert\prod^3_{i=1} \Big[\sum_s \lambda_s^{(i)} \Cal X_{\tau_s^{(i)}}\Big]\Big\Vert_{L^r}\leq c\delta^{\frac 3r} \prod^3_{i=1}
\Big[\sum_s|\lambda_s^{(i)}\Big].\tag 8.22
$$
Application of (8.22) to (8.19) with $r=\frac q6 >\frac 59$ implies, by (8.18)
$$
(8.7), (8.14), (8.19) < C\lambda^{3/2}\Big( \frac 1{\sqrt\lambda}\Big)^3 \prod^3_{i=1} \Vert f_i\Vert_2^{q/3}\tag 8.23
$$
and in view of the initial rescaling, (8.2) follows.

\bigskip
\noindent
{\bf Appendix: Upsilon Removal Lemmas}

We consider first the restriction (or extension) problem.

What follows is basically a modification of Theorem 1.2 in [T2] on deriving global restriction estimates from local ones. 
A significant difference is that instead of considering bounds of the type $(\gamma>0)$
$$
\Vert\hat f|_S\Vert_{L^{p} (d\sigma)} \lesssim R^\gamma\Vert f\Vert_{L^p(B_R)}\tag 1
$$
for $f\in L^p(\Bbb R^n), \supp f\subset B_R$, we start from a local inequality of the form
$$
\Vert\hat f|_S\Vert_{L^1(d\sigma)} \lesssim R^\gamma\Vert f\Vert_{L^p(B_R)}.\tag 2
$$
Compared with the argument from [T2], this will require additional involvement of the Maurey-Nikishin factorization theorem.

\proclaim
{Lemma A1} Assuming $1<p<2$,  $0<\gamma\ll 1$ and (2). Then
$$
\Vert\hat f|_S\Vert_{L^1(d\sigma)} \lesssim \Vert f\Vert_{p_1}\tag 3
$$
for $f\in L^{p_1}(\Bbb R^n)$ and
$$
\frac 1{p_1}>\frac 1p +\frac C{\log \frac 1\gamma}.\tag 4
$$
\endproclaim

In particular, if $(1)$ holds for arbitrary small $\gamma >0$, the global inequality (3) will hold for any $p_1< p$.

We start by dualizing (2), implying that the operator
$$
T:L^\infty (S, d\sigma)\to L^{p'}(B_R):\vp\to\widehat{\vp\sigma}|_{B_R}\qquad \Big(p' =\frac p{p-1}\Big)
$$
satisfies $\Vert T\Vert <R^\gamma$.
Hence, from the theory of absolutely summary operators, fixing any $r>p'>2$, there is a probability measure $\mu$ on $S$, such that
$$
\Vert\widehat{\vp\sigma}\Vert_{L^{p'}(B_R)} \lesssim R^\gamma\Vert\vp\Vert_{L^r(d\mu)}.\tag 5
$$
There is no harm to assume $\frac {d\mu}{d\sigma} >\frac 12$.

We first enforce some smoothness for the density.
Let $\tau:S\to S$ be any diffeomorphism that is $\frac 1R$-close to the identity.
Then, for $|x|<R$, a change of variables gives
$$
\align
\widehat{(\vp\circ \tau)\sigma} (x) &=\int\vp(\tau(\xi)\big) e(x.\xi)\sigma(d\xi)= \\
&\int \vp(\xi') e\big(x.\tau^{-1}(\xi')\big)\Delta (\xi') \sigma(dx') \  \text { where } |1-\Delta|\lesssim \frac 1R\\
&=\int(\Delta\vp)(\xi') e(x.\xi')\sigma(d\xi')\tag 6\\
&+O\big\{\sum_{j\geq 1} \frac 1{j!} \Big|\int(\Delta\psi_j\vp)(\xi') e(x\xi')\sigma(d\xi')\Big|\Big\}\tag 7
\endalign
$$
where (7) is obtained by Taylor expansion of $e\big(x.(\tau^{-1}(\xi')-\xi')\big)$ and \hfill\break
$|\psi_j(\xi')|<(R|\tau^{-1}(\xi')-\xi'|)^j< 1$ by 
assumption on $\tau$.
Hence
$$
|T(\vp\circ\tau)|\leq |T(\vp\Delta)|+\sum_{j\geq 1}\frac 1{j!} |T(\Delta\psi_j\vp)|
$$
and applying (5)
$$
\Vert T(\vp\circ\tau)\Vert_{L^{p'}(B_R)} \lesssim R^\gamma\Vert\vp\Vert_{L^r(d\mu)}.
$$
Replacing $\vp$ by $\vp\circ \tau^{-1}$, we obtain
$$
\Vert T\vp\Vert_{p'} \lesssim R^\gamma\Vert\vp\circ\tau^{-1}\Vert _{L^r(d\mu)} =R^\gamma\Vert\vp\Vert_{L^r(d\mu_\tau)}
$$
with $\mu_\tau =(\tau^{-1})_*[\mu]$.
Averaging over $\tau$ as above allows us to smoothen out $\mu$ at scale $\frac 1R$ and replace $\mu$ by a probability measure $\mu'$ on $S$, $\mu\ll\sigma$
and $\frac{d\mu'}{d\sigma}=\rho\geq \frac 12$ with $\rho$ smooth at scale $\frac 1R$.
Thus we have
$$
\Vert\widehat{\vp\sigma}\Vert_{L^{^{p'}}(B_R)} \leq R^\gamma\Vert\vp\rho^{1/r}\Vert_{L^r(d\sigma)}\tag 8
$$
and dualizing
$$
\Vert\hat f\rho^{-1/r}\big|_S\Vert_{L^{r'}(d\sigma)} \leq R^\gamma\Vert f\Vert_p \text { if } \supp f\subset B_R.\tag 8$'$
$$
In what preceeds, we fixed $R\geq 1$.
Note that $\rho$ depends on $R$.

Following [T2], define a finite collection of balls $\{B(a_\alpha, R)\}^N_{\alpha =1}$ in $\Bbb R^n$ as `sparse' if for $\alpha\not=\alpha'$
$$
|a_\alpha -a_{\alpha'}|> (NR)^C\tag 9
$$
($C$ some constant to specify).

Let now $\supp f\subset\bigcup_\alpha B(a_\alpha, R)$, i.e.
$$
f=\sum^N_{\alpha=1} f_\alpha(x - a_\alpha) \text { with } \supp f_\alpha \subset B_R.
$$
Hence
$$
\hat f(\xi)=\sum e(a_\alpha.\xi) \hat f_\alpha (\xi)
$$
and since $\Vert\vp\Vert_1=1$
$$
\Vert\hat f\big|_S\Vert_{L^1(d\sigma)}\leq\Big\Vert\Big[\sum e(a_\alpha.\xi)\hat f_\alpha (\xi)\Big] \rho^{-\frac 1r}(\xi)\Big\Vert
_{L^{r'}(d\sigma)}.\tag 10
$$
Note that by our construction of $\rho$, the function $ g_\alpha =\hat f_\alpha. \rho^{-\frac 1r}\big|_S$ is smooth at scale $\frac 1R$.
The sparsity of $\{a_\alpha\}$ allows then to estimate
$$
\Big\Vert \sum^N_{1} e(a_\alpha.\xi) g_\alpha(\xi)\Big\Vert_{L^{r'}(d\sigma)} \leq  2\Big(\sum\Vert g_\alpha\Vert^{r'}_{L^{r'}
(d\sigma)}\Big)^{1/{r'}}.\tag 11
$$
This is basically Lemma 3.2 in [T2] and we include the argument for completeness sake.

Establish (11) by interpolation.

More precisely, the claim will follow from an inequality for $1\leq s\leq 2$
$$
\Big\Vert\sum^N_1 e(a_\alpha. \xi) (\tilde\vp_\alpha *P_{\frac 1R}) (\xi)\Big|_S \Big\Vert _{L^s(d\sigma)} \lesssim \Big(\sum\Vert
\vp_\alpha\Vert^s_{L^s(d\sigma)}\Big)^{\frac 1s}\tag 12
$$
where $\{\vp_\alpha\}$ are arbitrary functions in $L^s(S, d\sigma)$, $^\sim$ denotes a well-behaved extension operator from $L^*(S)\to L^*(\Bbb
R^n)$ (take for instance the harmonic extension) and $P_{\frac 1R}$ is an $\frac 1R$-approximate identity.

For $s=1$, (12) is trivial from triangle inequality and since $\Vert(\tilde\vp * P_{\frac 1R})\big|_S \Vert_1\lesssim \Vert\vp\Vert_1$.

For $s=2$, we obtain for the square of the left side of (12)
$$
\sum_1^N\Vert\vp_\alpha\Vert^2_2 +\sum_{\alpha\not=\alpha'} \Big|\int e\big((a_\alpha -a_{\alpha'}).\xi\big) 
(\tilde\vp_\alpha * P_{\frac 1R}) (\xi) . \overline{(\tilde\vp_{\alpha'} * P_{\frac 1R})} (\xi)\sigma(d\xi)\Big|\tag 13
$$
and show that the contribution of the off-diagonal is small.

Denoting $\Phi_\alpha =\tilde\vp_\alpha * P_{\frac 1R}$, we may assume $\supp \hat\Phi_\alpha \subset B_R$ so that clearly, invoking the decay of $\hat\sigma$
and the fact that $|a_\alpha -a_{\alpha'}|\gg R$
$$
\align
&\Big|\int e\big((a_\alpha -a_{\alpha'}).\xi\big) \Phi_\alpha (\xi) \overline{\Phi_{\alpha'}}(\xi)\sigma(d\xi)\Big|\lesssim\\
&{}\\
&\frac 1{|a_\alpha - a_{\alpha'}|^{\frac{n-1}2}} \ \Vert\hat\phi_\alpha\Vert_1 \ \Vert\hat\phi_{\alpha'}\Vert_1\lesssim \frac
{R^n}{|a_\alpha-a_{\alpha'}|^{\frac {n-1}2}} \Vert\phi_\alpha\Vert_2 \ \Vert\phi_{\alpha'}\Vert_2\\
&{}\\
&\lesssim \frac {R^n}{|a_\alpha -a_{\alpha'}|^{\frac {n-1} 2}} \ \Vert \vp_\alpha\Vert_2 \Vert\vp_{\alpha'}\Vert_2.
\endalign
$$
Consequently, the second term in (13) is bounded by the first, provided
$$
\max_\alpha\sum_{\alpha'\not= \alpha} \frac 1{|a_\alpha -a_{\alpha'}|^{\frac{n-1}2}} < R^{-n}.
$$

This will be ensured if we require for $\alpha\not= \alpha'$
$$
|a_\alpha -a_{\alpha'}|> N^{\frac {n+1}{n(n-1)}} R^{\frac {2n}{n-1}}\tag 14
$$
as implied by (9) for $C$ large enough.
Then (12) will hold for $s=2$ and hence for $1\leq s\leq 2$.
Thus we proved (11).

Application of (11) with $g_\alpha =\hat f_\alpha.\rho^{-\frac 1r}|_S$ and invoking (8$'$) implies that
$$
\Vert
\hat f|_S\Vert_{L^1(d\sigma)} \lesssim R^\gamma\Big(\sum^N_{\alpha=1} \Vert f_\alpha\Vert_p^{r'}\Big)^{\frac 1{r'}}\lesssim
R^\gamma N^{\frac 1{r'}-\frac 1p}\Vert f\Vert_p\tag 15
$$
(recall that $r>p'$ is arbitrary).

Thus inequality (15) holds provided supp $f$ is contained in a sparse collection of $N$ balls of radius $R$.

The next ingredient is the following covering lemma (Lemma 3.3) from [T2].

\proclaim
{Lemma A2}
Suppose $E\subset \Bbb R^n$ is a finite union of 1-cubes and take $0<\delta<1$.
Then there exist $O\big(\frac 1\delta|E|^\delta)$ sparse collections of balls that cover $E$, such that the balls in each collection have radius
at most $O\big( \, |E|^{C^{1/\delta}}\big)$.
\endproclaim

Of course the number of balls in each collection is trivially bounded by $|E|$.

Assume $\supp f\subset E$ and apply Lemma A2 to $E$ (assumed a union of 1-cubes).
Hence
$$
E\subset\bigcup_{j\lesssim\frac 1\delta|E|^\delta} \  \bigcup_{a\in\Cal E_j} B(a, R_j)
$$
with $R_j\lesssim |E|^{C^{1/\delta}}$ and $\{B(a, R_j); a\in\Cal E_j\}$ sparse for each $j;\# \Cal E_j \leq N=|E|$.

Writing $f=\sum f_j, f_j= f\big|_{\bigcup_{a\in\Cal E_j}B(a, R_j)}$, application of inequality (15) to each $f_j$ implies
$$
\Vert \hat f\big|_S \Vert_{L^1(d\sigma)} \lesssim \frac 1\delta|E|^{\gamma C^{1/\delta}+\delta} N^{\frac 1{r'}-\frac 1p} \Vert f\Vert_p.\tag 16
$$
Taking $\delta\sim\frac 1{\log\frac 1\gamma}$ and $r< p'+\frac 1{\log \frac 1\gamma}$, we conclude that
$$
\Vert\hat f\big|_S\Vert _{L^1(d\sigma)} \lesssim_\gamma |E| ^{\frac C{\log\frac 1\gamma}} \Vert f\Vert_p.\tag 17
$$
Let $p_1<p$ and $f\in L^{p_1}(\Bbb R^n)$, $\Vert f\Vert_{p_1}\leq 1$, which we assume constant on $c$-cubes $(c\sim 1)$.

Decompose in level sets
$$
f=\sum_{k\geq 0} f\big|_{[2^{-k-1} \leq |f|< 2^{-k}]} =\sum f_k
$$
with $\supp f_k =E_k, E_k$ a union of $N_k$ $c$-cubes and $2^{-kp_1} N_k \lesssim 1$.

>From (17)
$$
\Vert\hat f_k\big|_S \Vert_{L^1(d\sigma)}\lesssim N_k^{\frac c{\log\frac 1\gamma}} \Vert f_k\Vert_p\lesssim 2^{k[\frac {cp_1}{\log \frac 1\gamma}+ \frac {p_1}
p - 1]}
$$
and therefore
$$
\Vert\hat f\big|_S\Vert_{L^1(d\sigma)} < C_\gamma\tag 18
$$
provided
$$
\frac C{\log\frac 1\gamma} < 1-\frac{p_1} p\tag 19
$$
which amounts to condition (4).

Arguing like in [T2], we showed that (18) holds for any function $f\in L^{p_1} (\Bbb R^n)$ of the form $f=\sum_{\xi\in\Cal L} \lambda_\xi
1_{B(\xi, c)}$ with $\sum|\lambda_\xi|^{p_1}\leq 1$ and $\Cal L$ a (shifted) 1-lattice.
Taking $c>0$ a sufficiently small constant as to ensure that $\widehat{1_{B_{(0, c)}}}$ is positive on $S$, it follows that
$$
\Big\Vert\Big[\sum_{\xi\in\Cal L} \lambda_\xi e(x.\xi)\Big]\Big|_S\Big\Vert_{L^1(d\sigma)}< C\Big(\sum|\lambda_\xi|^{p_1}\Big)^{\frac 1{p_1}}.\tag 20
$$
Another averaging over translates $\Cal L$ of the $\Bbb Z^n$-lattice gives (3).

This completes the proof of Lemma A1.

\bigskip

Next, we prove the upsilon-removal lemma in the variable coefficient multilinear case. Recall the setting.

Consider $T_\lambda$ and in (1.4), (1.5) with fixed, large $\lambda$ and define
$$
(Tf)(x) =\int e^{i\phi(x, y)} f(y)dy\tag 21
$$
with
$$
\phi(x, y) =x_1 y_1+\cdots +x_{n-1} y_{n-1} +x_n \big(\langle Ay, y\rangle +O(|y|^3)\big)
+\lambda\phi_\nu\Big(\frac xA, y\Big)\tag 22
$$
as in \S5, where $|x|=o(\lambda), |y|=o(1)$ and $A$ non-degenerate.

Let $2\leq k\leq n$ and $U_1, \ldots, U_k$ fixed balls in $y$-space satisfying the transversality condition (5.6).
For $j= 1, \ldots k$, denote
$$
T_jf=\int_{U_j}\ e^{i\phi(x, y)}f(y)dy.\tag 23
$$
Clearly the [BCT] result implies that if $1< R< o(\lambda)$, then
$$
\Big\Vert\Big(\prod^k_{j=1} |T_jf_j|\Big)^{\frac 1k}\Big\Vert_{L^q(B_R)} \ll R^\ve \Big(\prod^k_{j=1} \Vert f_j\Vert_2\Big)^{1/k}\tag 24
$$
with $q =\frac{2k}{k-1}$ and $B_R=B(0, R)$.
This statement is also easily seen to imply (24) with $B_R=B(a, R)$ any $R$-ball with $|a|=o(\lambda)$.

Our aim is to remove the $R^\ve$-factor at the cost of increasing slightly the exponent $q$.
Thus

\proclaim
{Lemma A3}
Under the above assumptions and taking $q_1>\frac {2k}{k-1}$, we have an inequality
$$
\Big\Vert\Big(\prod^k_{j=1} |T_jf_j|\big)^{\frac 1k}\Big\Vert_{q_1} \leq C_{n, k, q_1}\Big(\prod^k_{j=1} \Vert f_j\Vert_2\Big)^{1/k}.\tag 25
$$
\endproclaim
\noindent
\big(Note that we do not claim removal of the $\lambda^\ve$-factor in Theorem 6.2 from [BCT], as the context of our Lemma A3 is more restrictive,
since the $T_j$-operators are given by (22), (23)\big.)

Let $\Vert f_j\Vert_2 =1$ and $F=\big(\prod^k_{1} |T_jf_j|\big)^{\frac 1k}$.

Let $E\subset\Bbb R^d$ be obtained as  union of a sparse collection of $R$-balls $B(a_\alpha, R)$, $|a_\alpha|=o(\lambda)$
with $\alpha =1,\ldots, N$.
We will show that
$$
\Vert F|_E\Vert_q < C_\ve R^\ve.\tag 26
$$
Using Lemma A2, this will imply that for $E'\subset\Bbb R^n$ any finite union of 1-cubes we have
$$
\Vert F\big|_{E'}\Vert_q< \frac 1\delta C_\ve |E'|^{\delta+\ve C^{1/\delta}}\tag 27
$$
with $\delta>0$ a parameter. Hence, for all $\ve<0$
$$
\Vert F\big|_{E'}\Vert_q< C_\ve' |E'|^\ve\tag 28
$$
from where one easily deduces that $\Vert F\Vert_{q_1}< C_{q_1}$ for $q_1>q$.

Let $E=\bigcup B(a_\alpha, R)$ be as above and fix $\alpha$.
Write for $x\in B(a_\alpha, R)$
$$
(Tf)(x) =\int e^{i[\phi(x, y)-\phi(a_\alpha, y)]} \big (e^{i\phi(a_\alpha, y)} f(y)\big) \ \omega_j(y)dy\tag 29
$$
with $\omega_j$ a smooth localization on $U_j$.

Denoting $g(y) = e^{i\phi(a_\alpha, y)} f(y)$,
$$
(29) =\int\Big[\int e^{i[\phi(x, y)-\phi(a_\alpha, y)+\xi y]} \omega_j(y)dy\Big] \hat g(\xi) d\xi.\tag 30
$$

Since $|\nabla_y[\phi(x, y)-\phi(a_{\alpha}, y)]|\lesssim |x-a_\alpha|\lesssim R$, we may clearly replace in (30) the function $g$ by
$P_{R_1} g=(\hat g \eta_{_{R_1}})^\vee$, denoting $\eta_{_{R_1}} (z) =\eta (\frac z{R_1})$ where $0\leq\eta\leq 1$ is a smooth bumpfunction with
$\eta(0)=1$, and taking say
$$
R_1= 100NR.\tag 31
$$

The remaining contribution to (30) will then indeed by $L^\infty$-bounded by $0\big((NR)^{-C}\big)$.

Defining
$$
f_\alpha = e^{-i\phi(a_\alpha, y)} P_{R_1} \big( e^{i\phi(a_\alpha, y)} f\big)
$$
we can thus replace $Tf$ by $Tf_\alpha$ on $B(a_\alpha, R)$.
Note that $|f_{j, \alpha}|\leq |f_j| * |\overset\vee \to\eta_{R_1}|$ may clearly be assumed supported by $U_j$.

Estimate
$$
\align
\Vert F\big|_E\Vert^q_q&=\sum_\alpha\Vert F\big|_{B(a_\alpha, R)}\Vert_q^q\\
&=\sum_\alpha\Big\Vert\Big(\prod_j\big| T_j(f_{j, \alpha})|\Big)^{\frac 1k} \Big\Vert^q_{L^q(B(a_\alpha, R))}+o(1)\\
&\overset 24\to\leq C_\ve R^{q\ve}\sum_\alpha\Big(\prod_j \Vert f_{j, \alpha}\Vert_2\Big)^{q/k}+o(1)\\
&< C_\ve R^{q\ve} \max_j\Big[\sum_\alpha \Vert f_{j, \alpha} \Vert_2^q\Big]+ o(1).\tag 32
\endalign
$$
Since $q>2$,
$$
\Big(\sum_\alpha\Vert f_\alpha\Vert_1^q\Big)^{\frac 1q} \leq\Big(\sum\Vert f_\alpha\Vert^2_2\Big)^{\frac 12} =\Big(\sum_\alpha
\Vert P_{R_1} (e^{i\phi(a_{\alpha, y})} f)\Vert^2_2
\Big)^{1/2}.\tag 33
$$
To bound (33), take functions $\{\zeta_\alpha\}$ such that supp $\hat\zeta_\alpha \subset B(0, R_1)$ and $\sum_\alpha\Vert\zeta_\alpha\Vert^2_2=1$ 
and evaluate
$$
\sum_\alpha\langle e^{i\phi(a_{\alpha}, \cdot)} f, \zeta_\alpha\rangle \leq \Big\Vert\sum_\alpha e^{i\phi(a_\alpha, y)} \zeta_\alpha(y)\Big\Vert_2 \ \Vert
f\Vert_2.\tag 34
$$

For the off-diagonal terms $\alpha\not= \beta$
$$
|\langle e^{i\phi(a_{\alpha}, \cdot)} \zeta_{\alpha},  e^{i\phi(a_{\beta}, \cdot)} \zeta_\beta\rangle| =\Big|\int e^{i[\phi(a_\alpha, y)-\phi(a_\beta, y)]}
(\zeta_\alpha\bar\zeta_\beta)(y)dy\Big|.\tag 35
$$
where
$$
\align
\phi(a, y)-\phi(a', y) &=(a_1-a_1') y_1+\cdots+ (a_{d-1}-a_{d-1}') y_{d-1} +(a_d-a_{d'})\big(\langle Ay, y\rangle +O({|y|^3})\big)\\
&+\lambda \Big[\phi_\nu\Big(\frac
a\lambda, y\Big) -\phi_\nu \Big(\frac {a'}\lambda, y\Big)\Big]
\endalign
$$
satisfies either
$$
|\nabla_y [\phi(a, y)-\phi(a', y)]|\gtrsim |a-a'|
$$
or
$$
|\det D_y^2 [\phi(a, y)-\phi(a', y)]|\gtrsim |a-a'|^{n-1}.
$$
Hence, recalling the sparsity assumption $|a_\alpha-a_\beta|>(NR)^C\gg R_1$, it follows that
$$
(34) \lesssim |a_\alpha-a_\beta|^{-\frac{n-1}2} \Vert\hat\zeta_\alpha\Vert_1 \  \Vert\hat\zeta_\beta\Vert_1 \lesssim
R_1^{n-1} (NR)^{-C} \Vert\zeta_\alpha\Vert_2 \ \Vert\zeta_\beta\Vert_2.\tag 36
$$
Therefore (34) $\leq 2 \big(\sum \Vert\zeta_\alpha\Vert^2_2\big)^{1/2} \leq 2$ and (33) is bounded.
Inequality (26) now follows from (32), completing the proof of Lemma A3.

\enddocument